\documentclass[letterpaper,english]{amsart}
\usepackage{pslatex}
\usepackage[T1]{fontenc}
\usepackage[latin1]{inputenc}
\usepackage{fancyhdr}
\usepackage{graphics}
\makeatletter

\providecommand{\LyX}{L\kern-.1667em\lower.25em\hbox{Y}\kern-.125emX\@}

 \theoremstyle{plain}    
 \newtheorem{thm}{Theorem}[section]
 \numberwithin{equation}{section} 
 \numberwithin{figure}{section} 
 \theoremstyle{plain}    
 \newtheorem{cor}[thm]{Corollary} 
 \theoremstyle{plain}    
 \newtheorem{lem}[thm]{Lemma} 
 \theoremstyle{plain}    
 \newtheorem{prop}[thm]{Proposition} 
 \theoremstyle{definition}
 \newtheorem{defn}[thm]{Definition}
 \theoremstyle{definition}
  \newtheorem{example}[thm]{Example}
 \theoremstyle{remark}
 \newtheorem{rem}[thm]{Remark}
 \theoremstyle{remark}    
 \newtheorem{claim}[thm]{Claim}
 \theoremstyle{remark}    
 \newtheorem*{claim*}{Claim}
 \newcommand{\lyxaddress}[1]{
   \par {\raggedright #1 
   \vspace{1.4em}
   \noindent\par}}

\usepackage{verbatim} 
\usepackage{varioref}
\usepackage{amssymb}
\usepackage[mathscr]{eucal}
\usepackage{amsthm}
\usepackage[all]{xy}
\usepackage{hyperref}
\usepackage{rotating}
\usepackage{psfrag}
\usepackage{float}

\mathcode`\:="603A 

\DeclareSymbolFont{rsfs}{U}{rsfs}{m}{n}
\DeclareSymbolFontAlphabet{\mathrf}{rsfs}

\newcommand{\tmap}[1]{{\mathrm T}_{#1}}

\newcommand{\rank}{\operatorname{rank}}

\newcommand{\im}{\operatorname{im}}

\newcommand{\nth}[1]{$#1^{\mathrm{th}}$}
\newcommand{\tunder}[2]{\tmap{\underbrace{\scriptstyle#1}_{{\text{#2}}}}}

\newcommand{\tunderi}[2]{\tunder{1,\dots,#1,\dots,1}{\nth #2\ position}}

\newdir{ >}{{}*!/-5pt/@{>}}
\SelectTips{cm}{}

\makeatother

\begin{document}

\title{Cofree coalgebras over operads}

\author{Justin R. Smith}



\lyxaddress{Department of Mathematics and Computer Science\\
Drexel University\\
Philadelphia,~PA 19104}

\email{jsmith@mcs.drexel.edu}

\urladdr{http://vorpal.mcs.drexel.edu}

\date{\today}

\begin{abstract}
This paper explicitly constructs cofree coalgebras over operads in
the category of DG-modules. Special cases are also considered in which
the general expression simplifies (such as the pointed, irreducible
case). 
\end{abstract}
\maketitle
\tableofcontents{}

\newcommand{\ring}{R}

\newcommand{\integers}{\mathbb{Z}}

\newcommand{\betabar}{\bar{\beta }}
 
\newcommand{\desusp}{\downarrow }

\newcommand{\susp}{\uparrow }

\newcommand{\cobar}{\mathcal{F}}

\newcommand{\bigp}{\mathop {{\prod }^{\prime }}}

\newcommand{\mfrac}{\mathfrak{M}}

\newcommand{\coend}{\mathrm{CoEnd}}

\newcommand{\ainfty}{A_{\infty }}

\newcommand{\coassoc}{\mathrm{Coassoc}}

\newcommand{\trm}{\mathrm{T}}

\newcommand{\tfr}{\mathfrak{T}}

\newcommand{\tabbr}{\hat{\trm }}

\newcommand{\Tabbr}{\hat{\tfr }}

\newcommand{\afr}{\mathfrak{A}}

\newcommand{\homz}{\mathrm{Hom}_{\ring }}

\newcommand{\homa}{\mathrm{Hom}}

\newcommand{\zend}{\mathrm{End}}

\newcommand{\rs}[1]{\mathrm{R}S_{#1 }}

\newcommand{\highprod}[1]{\bar{\mu }_{#1 }}

\newcommand{\slength}[1]{|#1 |}

\newcommand{\barcs}{\bar{\mathcal{B}}}

\newcommand{\ubarcs}{\mathcal{B}}

\newcommand{\zs}[1]{\mathbb{Z}S_{#1 }}

\newcommand{\homzs}[1]{\mathrm{Hom}_{\ring S_{#1 }}}

\newcommand{\zpi}{\mathbb{Z}\pi }

\newcommand{\D}{\mathfrak{D}}

\newcommand{\ahat}{\hat{\mathfrak{A}}}

\newcommand{\cbar}{{\bar{C}}}

\newcommand{\cf}[1]{\mathscr{C}(#1 )}

\newcommand{\ddelta}{\dot{\Delta }}

\newcommand{\dimlimiter}{\triangleright }

\section{Introduction}

\newdir{ >}{{}*!/-5pt/@{>}}

We begin with the definition of the object of this paper:

\begin{defn}
\label{def:cofreecoalgebra}Let \( R \) be a commutative ring with
unit and \( C \) be an \( R \)-module. Then a coalgebra \( G \)
will be called \emph{the cofree coalgebra} \emph{cogenerated by} \( C \)
if
\begin{enumerate}
\item there exists a morphism of \( R \)-modules\[
\varepsilon :G\to C\]
called the \emph{cogeneration map},
\item given any coalgebra \( D \) and any morphism of \( R \)-modules\[
f:D\to C\]
 there exists a \emph{unique} morphism of coalgebras\[
\hat{f}:D\to G\]
 that makes the diagram \[\xymatrix{{D} \ar[r]^-{\hat{f}}\ar[rd]_{f}& {G}\ar[d]^{\varepsilon} \\ & {C}}\] commute.
\end{enumerate}
If \( \mathscr{V} \) is an operad (see definition~\ref{def:operad})
and \( C \) is a \( \ring  \)-free DG-module, then the same definition
holds for coalgebras and coalgebra-morphisms (see definition~\ref{dia:coalgdef})
over \( \mathscr{V} \).
\end{defn}
\begin{rem}
If they exist, it is not hard to see that cofree coalgebras must be
\emph{unique} up to an isomorphism.
\end{rem}
Constructions of free \emph{algebras} satisfying various conditions
(associativity, etc.) have been known for many years: One forms a
general algebraic structure implementing a suitable {}``product''
and forms the quotient by a sub-object representing the conditions.
Then one shows that these free algebras map to any other algebra satisfying
the conditions. For instance, it is well-known how to construct the
free \emph{algebra} over an operad --- see \cite{Kriz-May}. 

The construction of cofree coalgebras is dual to this, although Thomas
Fox showed (see \cite{T.Fox:1993, T.Fox:2001}) that they are considerably
more complex than free algebras. Definition~\ref{def:cofreecoalgebra}
implies that a cofree coalgebra cogenerated by a \( R \)-module,
\( C \), must contain isomorphic images of \emph{all possible} coalgebras
over \( C \). 

Operads (in the category of graded groups) can be regarded as {}``systems
of indices'' for parametrizing operations. They provide a uniform
framework for describing many classes of algebraic objects, from associative
algebras and coalgebras to Lie algebras and coalgebras. 

In recent years, there have been applications of operads to quantum
mechanics and homotopy theory. For instance, Steenrod operations on
the chain-complex of a space can be codified by making this chain-complex
a coalgebra over a suitable operad.

The definitive references on cofree coalgebras are the the book \cite{Sweedler:1969}
and two papers of Fox. Sweedler approached cofree coalgebras as a
kind of dual of free algebras, while Fox studied them \emph{ab initio,}
under the most general possible conditions. 

In \S~\ref{sec:generalconstruction}, we describe the cofree coalgebra
over an operad and prove that it has the required properties. Theorem~\ref{th:cofreegeneral}
gives our result. 

In \S~\ref{sec:special} we consider consider special cases such as
the pointed irreducible case in which the coproduct is dual to the
operad compositions --- see \ref{th:pointedirreduciblecofree} and
\ref{prop:pointedirreducnonnegative}.

\section{Operads}

\subsection{Notation and conventions}

Throughout this paper, \( R \) will denote a commutative ring with
unit. All tensor-products will be over \( R \) so that \( \otimes =\otimes _{R} \).

\begin{defn}
\label{def:degmap} Let \( C \) and \( D \) be graded \( \ring  \)-modules.
A map of graded modules \( f:C_{i}\to D_{i+j} \) will be said to
be of degree \( j \). 
\end{defn}
\begin{rem}
For instance the \textit{differential} of a DG-module will be regarded
as a degree \( -1 \) map. 
\end{rem}
We will make extensive use of the Koszul Convention (see~\cite{Gugenheim:1960})
regarding signs in homological calculations:

\begin{defn}
\label{def:basic} \label{def:koszul} If \( f:C_{1}\to D_{1} \),
\( g:C_{2}\to D_{2} \) are maps, and \( a\otimes b\in C_{1}\otimes C_{2} \)
(where \( a \) is a homogeneous element), then \( (f\otimes g)(a\otimes b) \)
is defined to be \( (-1)^{\deg (g)\cdot \deg (a)}f(a)\otimes g(b) \). 
\end{defn}
\begin{rem}
This convention simplifies many of the common expressions that occur
in homological algebra --- in particular it eliminates complicated
signs that occur in these expressions. For instance the differential,
\( \partial _{\otimes } \), of the tensor product \( C\otimes D \)
is \( \partial _{C}\otimes 1+1\otimes \partial _{D} \).

If \( f_{i} \), \( g_{i} \) are maps, it isn't hard to verify that
the Koszul convention implies that \( (f_{1}\otimes g_{1})\circ (f_{2}\otimes g_{2})=(-1)^{\deg (f_{2})\cdot \deg (g_{1})}(f_{1}\circ f_{2}\otimes g_{1}\circ g_{2}) \).
\end{rem}
Another convention that we will follow extensively is tensor products,
direct products, etc. are of \emph{graded modules.} 

\emph{Powers} of DG-modules, such as \( C^{n} \) will be regarded
as iterated \( \ring  \)-tensor products: \[
C^{n}=\underbrace{C\otimes \cdots \otimes C}_{n\, \mathrm{factors}}\]

\subsection{Definitions}

Before we can define operads, we need the following:

\begin{defn}
\label{def:tmap}Let \( \sigma \in S_{n} \) be an element of the
symmetric group and let \( \{k_{1},\dots ,k_{n}\} \) be \( n \)
nonnegative integers with \( K=\sum _{i=1}^{n}k_{i} \). Then \( T_{k_{1},\dots ,k_{n}}(\sigma ) \)
is defined to be the element \( \tau \in S_{K} \) that permutes the
\( n \) blocks \[
(1,\dots ,k_{1}),(k_{1}+1,\dots ,k_{1}+k_{2})\dots (K-K_{n-1},\dots ,K)\]
as \( \sigma  \) permutes the set \( \{1,\dots ,n\} \).
\end{defn}
\begin{rem}
Note that it is possible for one of the \( k \)'s to be \( 0 \),
in which case the corresponding block is empty. 
\end{rem}
The standard definition (see \cite{Kriz-May}) of an operad in the
category of DG-modules is:

\begin{defn}
\label{def:operad} A sequence of differential graded \( \ring  \)-free
modules, \( \{\mathscr{V}_{i}\} \), will be said to form a \emph{DG-operad}
if they satisfy the following conditions:
\begin{enumerate}
\item there exists a \emph{unit map} (defined by the commutative diagrams
below) \[
\eta :\ring \to \mathscr{V}_{1}\]

\item for all \( i>1 \), \( \mathscr{V}_{i} \) is equipped with a left
action of \( S_{i} \), the symmetric group. 
\item for all \( k\ge 1 \), and \( i_{s}\ge 0 \) there are maps \[
\gamma :\mathscr{V}_{k}\otimes \mathscr{V}_{i_{1}}\otimes \cdots \otimes \mathscr{V}_{i_{k}}\otimes \to \mathscr{V}_{i}\]
 where \( i=\sum _{j=1}^{k}i_{j} \).

The \( \gamma  \)-maps must satisfy the following conditions:

\end{enumerate}
\begin{description}
\item [Associativity]the following diagrams commute, where \( \sum j_{t}=j \),
\( \sum i_{s}=i \), and \( g_{\alpha }=\sum _{\ell =1}^{\alpha }j_{\ell } \)
and \( h_{s}=\sum _{\beta =g_{s-1}+1}^{g_{s}}i_{\beta } \):  \begin{equation}\xymatrix@C+20pt{ {\mathscr{V}_{k}\otimes\left(\bigotimes_{s=1}^{k}\mathscr{V}_{j_{s}}\right) \otimes \left( \bigotimes_{t=1}^{j}\mathscr{V}_{i_{t}} \right) } \ar[r]^-{\gamma\otimes\mathrm{Id}} \ar[dd]_{\mathrm{shuffle}} & {\mathscr{V}_{j}\otimes\left(\bigotimes_{t=1}^{j}\mathscr{V}_{i_{t}}\right)} \ar[d]^{\gamma} \\
&{\mathscr{V}_{i}} \\ {\mathscr{V}_{k}\otimes\left(\bigotimes_{t=1}^{k}\mathscr{V}_{j_{t}}\otimes\left(\bigotimes_{q=1}^{j_{t}}\mathscr{V}_{i_{g_{t-1}+q}}\right)\right) }\ar[r]_-{\mathrm{Id}\otimes(\otimes_{t}\gamma)}&{ \mathscr{V}_{k}\otimes\left(\bigotimes_{t=1}^{k}\mathscr{V}_{h_{t}}\right)}\ar[u]_{\gamma} }\label{dia:operadassociativity}\end{equation}

\item [Units]the following diagrams commute:  \begin{equation}\begin{array}{cc}\xymatrix{{\mathscr{V}_{k}\otimes{\integers}^{k}}\ar[r]^{\cong}\ar[d]_{\text{Id}\otimes{\eta}^{k}}&{\mathscr{V}_{k}}\\
{\mathscr{V}_{k}\otimes{\mathscr{V}_{1}}^{k}}\ar[ur]_{\gamma}&}& \xymatrix{{\integers\otimes\mathscr{V}_{k}}\ar[r]^{\cong}\ar[d]_{\eta\otimes\text{Id}}&{\mathscr{V}_{k}}\\
{\mathscr{V}_{1}\otimes\mathscr{V}_{k}}\ar[ur]_{\gamma}&}\end{array}\label{dia:operadunit}\end{equation} 
\item [Equivariance]the following diagrams commute: \begin{equation}\xymatrix@C+20pt{{\mathscr{V}_{k}\otimes\mathscr{V}_{j_{1}}\otimes\cdots\otimes\mathscr{V}_{j_{k}}}\ar[r]^-{\gamma}\ar[d]_{\sigma\otimes\sigma^{-1}}&{\mathscr{V}_{j}}\ar[d]^{\tmap{j_{1},\dots,j_{k}}(\sigma)}\\
{\mathscr{V}_{k}\otimes\mathscr{V}_{j_{\sigma(1)}}\otimes\cdots\otimes\mathscr{V}_{j_{\sigma(k)}}}\ar[r]_-{\gamma}&{\mathscr{V}_{j}}}\label{dia:operadequivariant}\end{equation}
where \( \sigma \in S_{k} \), and the \( \sigma ^{-1} \) on the
left permutes the factors \( \{\mathscr{V}_{j_{i}}\} \) and the \( \sigma  \)
on the right simply acts on \( \mathscr{V}_{k} \). See \ref{def:tmap}
for a definition of \( \tmap {j_{1},\dots ,j_{k}}(\sigma ) \). \[\xymatrix@C+20pt{{\mathscr{V}_{k}\otimes\mathscr{V}_{j_{1}}\otimes\cdots\otimes\mathscr{V}_{j_{k}}}\ar[r]^-{\gamma}\ar[d]_{\text{Id}\otimes\tau_{1}\otimes\cdots\tau_{k}}&{\mathscr{V}_{j}}\ar[d]^-{\tau_{1}\oplus\cdots\oplus\tau_{k}}\\
{\mathscr{V}_{k}\otimes\mathscr{V}_{j_{\sigma(1)}}\otimes\cdots\otimes\mathscr{V}_{j_{\sigma(k)}}}\ar[r]_-{\gamma}&{\mathscr{V}_{j}}}\] where
\( \tau _{s}\in S_{j_{s}} \) and \( \tau _{1}\oplus \cdots \oplus \tau _{k}\in S_{j} \)
is the block sum.
\end{description}
The individual \( \mathscr{V}_{n} \) that make up the operad \( \mathscr{V} \)
will be called its \emph{components.}
\end{defn}
For reasons that will become clear in the sequel, we follow the nonstandard
convention of using \emph{subscripts} to denote components of an operad
--- so \( \mathscr{V}=\{\mathscr{V}_{n}\} \) rather than \( \{\mathscr{V}(n)\} \).
Where there is any possibility of confusion with grading of a graded
groups, we will include a remark.

We will also use the term \emph{operad} for DG-operad throughout this
paper.

\begin{defn}
\label{def:unitaloperad}An operad, \( \mathscr{V} \), is called
\emph{unital} if \( \mathscr{V} \) has a \( 0 \)-component \( \mathscr{V}_{0}=\ring  \),
concentrated in dimension \( 0 \) and augmentations\[
\epsilon _{n}:\mathscr{V}_{n}\otimes \mathscr{V}_{0}\otimes \cdots \otimes \mathscr{V}_{0}=\mathscr{V}_{n}\to \mathscr{V}_{0}=\ring \]
 induced by their structure maps.
\end{defn}
\begin{rem}
The literature contains varying definitions of the terms discussed
here.

Our definition of unital and non-unital operad corresponds to that
in \cite{Kriz-May}. On the other hand, in \cite{Markl:1996} Markl
defines a \emph{unital} operad to have a \emph{unit} (i.e., the map
\( \eta :\ring \to \mathscr{V}_{1} \)) and calls operads meeting
the condition in definition~\ref{def:unitaloperad} \emph{augmented
unital.} None of Markl's operads have a \( 0 \)-component and his
definition of augmentation \emph{only} involves the \( 1 \)-component
(so that the {}``higher'' augmentation maps \( \epsilon _{n}:\mathscr{V}_{n}\to \ring  \)
do not have to exist).
\end{rem}

\subsection{The composition-representation}

Describing an operad via the \( \gamma  \)-maps and the diagrams
in \ref{def:operad} is known as the \( \gamma  \)-representation
of the operad. We will present another method for describing operads
more suited to the constructions to follow:

\begin{defn}
\label{def:ithcomposition}Let \( \mathscr{V} \) be an operad as
defined in \ref{def:operad}, let \( n,m \) be positive integers
and let \( 1\le i\le n \). Define\[
\circ _{i}:\mathscr{V}_{n}\otimes \mathscr{V}_{m}\to \mathscr{V}_{n+m-1}\]
\emph{the \( i^{\mathrm{th}} \) composition operation} on \( \mathscr{V} \),
to be the composite \[\xymatrix{ {\mathscr{V}_n \otimes \mathscr{V}_m} \ar@{=}[d]\\ {\mathscr{V}_n \otimes \ring^{i-1}\otimes\mathscr{V}_m\otimes\ring^{n-i}} \ar[d]^{1\otimes\eta^{i-1}\otimes 1\otimes \eta^{n-i}}\\ {\mathscr{V}_n \otimes {\mathscr{V}_1}^{i-1}\otimes\mathscr{V}_m\otimes{\mathscr{V}_1}^{n-i}}\ar[d]^{\gamma} \\ {\mathscr{V}_{n+m-1}}}\] 
\end{defn}
The \( \gamma  \)-maps defined in \ref{def:operad} and the composition-operations
uniquely determine each other. 

\begin{defn}
\label{def:lcomps}Let \( \mathscr{V} \) be an operad, let \( 1\le j\le n \),
and let \( \{\alpha _{1},\dots ,\alpha _{j}\} \) be positive integers.
Then define\[
L_{j}:\mathscr{V}_{n}\otimes \mathscr{V}_{\alpha _{1}}\otimes \cdots \otimes \mathscr{V}_{\alpha _{j}}\to \mathscr{V}_{n-j+\sum \alpha _{i}}\]
to be the composite \begin{equation}\xymatrix{{\mathscr{V}_n\otimes\mathscr{V}_{\alpha_1}\otimes\cdots\otimes\mathscr{V}_{\alpha_j}}\ar@{=}[d] \\ { \mathscr{V}_n\otimes( \mathscr{V}_{\alpha_1}\otimes\cdots\otimes\mathscr{V}_{\alpha_j}\otimes\ring\otimes\cdots\otimes\ring)}\ar[d]^{1 \otimes(1^j\otimes\eta^{n-j})}\\ {\mathscr{V}_n \otimes( \mathscr{V}_{\alpha_1}\otimes\cdots\otimes\mathscr{V}_{\alpha_j}\otimes\mathscr{V}_1\otimes\cdots\otimes\mathscr{V}_1)}\ar[d]^{\gamma}\\{\mathscr{V}_{n+\sum_{i=1}^j(\alpha_i-1)}} }\label{dia:itercomposition}\end{equation}
\end{defn}
\begin{rem}
Clearly, under the hypotheses above, \( L_{n}=\gamma  \).

Operads were originally called \emph{composition algebras} and defined
in terms of these operations --- see \cite{Gerstenhaber:1962}. 
\end{rem}
\begin{prop}
\label{prop:gammaequivcomps}Under the hypotheses of \ref{def:lcomps},
suppose \( j<n \). Then \begin{eqnarray*}
L_{j+1}=L_{j}\circ (\ast \circ _{j+1+\sum _{i=1}^{j}\alpha _{i}}\ast ): &  & \\
\mathscr{V}_{n}\otimes \mathscr{V}_{\alpha _{1}}\otimes \cdots \otimes \mathscr{V}_{\alpha _{j+1}} & \to  & \mathscr{V}_{n+\sum _{i=1}^{j+1}(\alpha _{i}-1)}
\end{eqnarray*}
In particular, the \( \gamma  \)-map can be expressed as an iterated
sequence of compositions and \( \gamma  \)-maps and the composition-operations
determine each other.
\end{prop}
\begin{rem}
We will find the compositions more useful than the \( \gamma  \)-maps
in studying algebraic properties of coalgebras over \( \mathscr{V} \). 

The map \( \gamma  \) and the composition-operations \( \{\circ _{i}\} \)
define the \( \gamma  \)- and the \emph{composition-representations}
of \( \mathscr{V} \), respectively. 
\end{rem}
\begin{proof}
This follows by induction on \( j \): it follows from the definition
of the \( \{\circ _{i}\} \) in the case where \( j=1 \). In the
general case, it follows by applying the associativity identities
and the identities involving the unit map, \( \eta :\ring \to \mathscr{V}_{1} \).
Consider the diagram \begin{equation}\xymatrix@R+10pt{{\mathscr{V}_n\otimes (\mathscr{V}_{\alpha_1}\otimes\cdots\otimes\mathscr{V}_{\alpha_j}\otimes\ring^{n-j}) \otimes\ring^{j+\sum_{i=1}^j(\alpha_i-1)}\otimes \mathscr{V}_{\alpha_{j+1}}\otimes\ring^{n-j-1} }\ar[d]|-{1 \otimes ( 1^{j}\otimes\eta^{n-j})\otimes\eta^{j+\sum_{i=1}^j(\alpha_i-1)} \otimes 1\otimes\eta^{n-j-1} }\\ {\mathscr{V}_n \otimes( \mathscr{V}_{\alpha_1}\otimes\cdots\otimes\mathscr{V}_{\alpha_t}\otimes\mathscr{V}^{n-j}_1) \otimes\mathscr{V}_1^{j+\sum_{i=1}^j(\alpha_i-1)}\otimes \mathscr{V}_{\alpha_{j+1}}\otimes\mathscr{V}_1^{n-j-1} }\ar[d]^{\gamma\otimes 1^{n+\sum_{i=1}^j(\alpha_i-1)}}\\{\mathscr{V}_{n+\sum_{i=1}^j(\alpha_i-1)}\otimes(\mathscr{V}_1^{j+\sum_{i=1}^j(\alpha_i-1)}\otimes \mathscr{V}_{\alpha_{j+1}}\otimes\mathscr{V}_1^{n-j-1}) }\ar[d]^{\gamma} \\ {\mathscr{V}_{n+\sum_{i=1}^{j+1}(\alpha_i-1)}} }\label{dia:itercomposition2}\end{equation}

The associativity condition implies that we can shuffle copies of
\( \mathscr{V}_{1} \) to the immediate left of the rightmost term,
and shuffle the \( \mathscr{V}_{1}\otimes \cdots \otimes \mathscr{V}_{1} \)
on the right to get a factor on the left of \[
\mathscr{V}_{\alpha _{1}}\otimes \mathscr{V}^{\alpha _{1}}_{1}\otimes \cdots \otimes \mathscr{V}_{\alpha _{j}}\otimes \mathscr{V}_{1}^{\alpha _{j}}\]
and one on the right of\[
\mathscr{V}_{1}\otimes \mathscr{V}_{\alpha _{j+1}}\]
(this factor of \( \mathscr{V}_{1} \) exists because \( j<n \))
and we can evaluate \( \gamma  \) on each of these before evaluating
\( \gamma  \) on their tensor product. The conclusion follows from
the fact that each copy of \( \mathscr{V}_{1} \) that appears in
the result has been composed with the unit map \( \eta  \) so the
left factor is\begin{eqnarray*}
\gamma (\mathscr{V}_{\alpha _{1}}\otimes \mathscr{V}^{\alpha _{1}}_{1})\otimes \cdots \otimes \gamma (\mathscr{V}_{\alpha _{j}}\otimes \mathscr{V}^{\alpha _{j}}_{1}) & = & \\
\mathscr{V}_{\alpha _{1}}\otimes \cdots \otimes \mathscr{V}_{\alpha _{j}} &  & 
\end{eqnarray*}
and the right factor is \[
\gamma (\mathscr{V}_{1}\otimes \mathscr{V}_{\alpha _{j+1}})=\mathscr{V}_{\alpha _{j+1}}\]
so the entire expression becomes\[
\gamma (\mathscr{V}_{n}\otimes \mathscr{V}_{\alpha _{1}}\otimes \cdots \otimes \mathscr{V}_{\alpha _{j+1}}\otimes \mathscr{V}_{1}^{n-j-1})\]
which is what we wanted to prove.
\end{proof}
The composition representation is complete when one notes that the
various diagrams in \emph{\ref{def:operad}} translate into the following
relations (whose proof is left as an exercise to the reader):

\begin{lem}
\label{def:compositionsassociativity}Compositions obey the identities

\( (a\circ _{i}b)\circ _{j}c=\left\{ \begin{array}{ll}
(-1)^{\dim b\cdot \dim c}(a\circ _{j-n+1}c)\circ _{i}b & \mathrm{if}\, i+n-1\le j\\
a\circ _{i}(b\circ _{j-i+1}c) & \mathrm{if}\, i\le j<i+n-1\\
(-1)^{\dim b\cdot \dim c}(a\circ _{j}c)\circ _{i+m-1}b & \mathrm{if}\, 1\le j<i
\end{array}\right.  \)

where \( \deg c=m \), \( \deg a=n \), and

\begin{equation}
\label{eq:operadcompequivar}
a\circ _{\sigma (i)}(\sigma \cdot b)=\tunderi {n}{i}(\sigma )\cdot (a\circ _{i}b)
\end{equation}

\end{lem}
Given compositions, we define \emph{generalized structure maps} of
operads.

\begin{defn}
\label{def:generalizedcomps}Let \( \mathscr{V} \) be an operad and
let \( \mathbf{u}=\{u_{1},\dots ,u_{m}\} \), be a list of symbols,
each of which is either a positive integer or the symbol \( \bullet  \).
We define the \emph{generalized composition} with respect to \( \mathbf{u} \),
denoted \( \gamma _{\mathbf{u}} \), by \[
\gamma _{\mathbf{u}}=\gamma \circ \bigotimes _{j=1}^{m}\iota _{j}:\mathscr{V}_{m}\otimes \mathscr{V}_{u_{1}}\otimes \cdots \otimes \mathscr{V}_{u_{m}}\to \mathscr{V}_{n}\]
where \[
n=\sum _{j=1}^{m}u_{j}\]
and we follow the convention that
\begin{enumerate}
\item \( \bullet =1 \) when used in a numeric context,
\item \( \mathscr{V}_{\bullet }=\ring  \),
\item \( \iota _{j}=\left\{ \begin{array}{ll}
1:\mathscr{V}_{u_{j}}\to \mathscr{V}_{u_{j}} & \mathrm{if}\, u_{j}\neq \bullet \\
\eta :\ring \to \mathscr{V}_{1} & \mathrm{otherwise}
\end{array}\right.  \)
\end{enumerate}
\end{defn}
\begin{rem}
If \( \{u_{k_{1}},\dots ,u_{k_{t}}\}\subset \{u_{1},\dots ,u_{m}\} \)
is the sublist of non-\( \bullet  \) elements, then \( \gamma _{\mathbf{u}} \)
is a map\[
\gamma _{\mathbf{u}}:\mathscr{V}_{m}\otimes \mathscr{V}_{k_{1}}\otimes \cdots \otimes \mathscr{V}_{k_{t}}\to \mathscr{V}_{n}\]

\end{rem}
\begin{lem}
\label{lem:generalizedcompassoc}Let \( \mathscr{V} \) be an operad,
let \( n,m,\alpha >0 \) let \( 1\le i\le n \) be integers, and let
\( \mathbf{u}=\{u_{1},\dots ,u_{n}\} \), \( \mathbf{v}=\{v_{1},\dots ,v_{m}\} \),
\( \mathbf{w}=\{w_{1},\dots ,w_{n+m-1}\} \) be lists of symbols as
in definition~\ref{def:generalizedcomps} with

\begin{eqnarray*}
u_{i} & \ne  & \bullet \\
u_{i} & = & \sum _{j=1}^{m}v_{j}\\
w_{j} & = & u_{j}\, \, \mathrm{if}\, \, j<i\\
w_{j} & = & v_{j-i+1}\, \, \mathrm{if}\, \, j\ge i\, \, \mathrm{and}\, \, j<i+m\\
w_{j} & = & u_{j-m+1}\, \, \mathrm{if}\, \, j\ge i+m\\
\alpha  & = & \sum _{j=1}^{n+m-1}w_{j}\\
 & = & \sum _{j=1}^{n}u_{j}
\end{eqnarray*}
Then the diagram \[\xymatrix{{\mathscr{V}_{n}\otimes\mathscr{V}_{m}\otimes\bigotimes_{k=1}^{n+m-1}\mathscr{V}_{w_k}}\ar[r]^{\circ_i\otimes1}\ar[d]_{1^{i-1}\otimes \gamma_{\mathbf{v}}\otimes1^{n-i}\circ s} & {\mathscr{V}_{n+m-1}\otimes\bigotimes_{k=1}^{n+m-1}\mathscr{V}_{w_k}}\ar[d]^{\gamma_{\mathbf{w}}} \\ {\mathscr{V}_n\otimes\bigotimes_{k=1}^n\mathscr{V}_{u_k}}\ar[r]_{\gamma_{\mathbf{u}}} &{\mathscr{V}_{\alpha}} }\] commutes,
where \( s \) is the shuffle map that sends \( \mathscr{V}_{m} \)
\( i-1 \) places to the right.
\end{lem}
\begin{rem}
The conditions on \( \mathbf{u} \), \( \mathbf{v} \), and \( \mathbf{w} \)
imply that \( \mathbf{w} \) is the result of replacing \( u_{i} \)
with the entire list \( \mathbf{v} \).
\end{rem}
\begin{proof}
This is a straightforward consequence of definition~\ref{def:generalizedcomps}
and the associativity condition in diagram~\ref{dia:operadassociativity}.
\end{proof}
Morphisms of operads are defined in the obvious way:

\begin{defn}
\label{def:operadmorphism} Given two operads \( \mathscr{V} \) and
\( \mathscr{W} \), a \emph{morphism} \[
f:\mathscr{V}\to \mathscr{W}\]
 is a sequence of chain-maps \[
f_{i}:\mathscr{V}_{i}\to \mathscr{W}_{i}\]
 commuting with all the diagrams in \ref{def:operad} or (equivalently)
preserving the composition operations in \ref{def:generalizedcomps}.
\end{defn}
Now we give some examples:

\begin{defn}
\label{def:mathfrakS0}The operad \( \mathfrak{S}_{0} \) is defined
via
\end{defn}
\begin{enumerate}
\item Its \( n^{\mathrm{th}} \) component is \( \zs{{n}} \) --- a chain-complex
concentrated in dimension \( 0 \). 
\item Its structure map is given by\[
\gamma (1_{S_{n}}\otimes 1_{S_{k_{1}}}\otimes \cdots \otimes 1_{S_{k_{n}}})=1_{S_{K}}\]
where \( 1_{S_{j}}\in S_{j} \) is the identity element and \( K=\sum _{j=1}^{n}k_{j} \).
This definition is extended to other values in the symmetric groups
via the equivariance conditions in \ref{def:operad}.
\end{enumerate}
\begin{rem}
This was denoted \( \mathrf {M} \) in \cite{Kriz-May}.
\end{rem}
Verification that this satisfies the required identities is left to
the reader as an exercise.

\begin{defn}
\label{def:sfrakfirstmention}Let \( \mathfrak{S} \) denote the operad
with components \( K(S_{n},1) \) --- the bar resolutions of \( \integers  \)
over \( \zs{{n}} \) for all \( n>0 \). See \cite{Smith:1994} for
formulas for the composition-operations.
\end{defn}
Now we define two important operads associated to any \( \ring  \)-module.

\begin{defn}
\label{def:coend} Let \( C \) be a DGA-module . Then the \emph{Coendomorphism}
operad, \( \coend (C) \), is defined to be the operad with component
of \( \rank i=\homz (C,C^{i}) \), with the differential induced by
that of \( C \) and \( C^{i} \). The dimension of an element of
\( \homz (C,C^{i}) \) (for some \( i \)) is defined to be its degree
as a map. If \( C \) is equipped with an augmentation\[
\varepsilon :C\to \ring \]
where \( \ring  \) is concentrated in dimension \( 0 \), then \( \coend (C) \)
is unital, with \( 0 \) component generated by \( \varepsilon  \)
(with the identification \( C^{0}=\ring  \)).
\end{defn}
\begin{rem}
One motivation for operads is that they model the iterated coproducts
that occur in \( \coend (\ast ) \). We will use operads as an algebraic
framework for defining other constructs that have topological applications.
\end{rem}

\subsection{Coalgebras over an operad}

\begin{defn}
\label{def:operadcomodule} Let \( \mathscr{V} \) be an operad and
let \( C \) be a DG-module equipped with a morphism (of operads)
\[
f:\mathscr{V}\to \coend (C)\]
 Then \( C \) is called a \emph{coalgebra} over \( \mathscr{V} \)
\emph{with structure map} \( f \).
\end{defn}
\begin{rem}
A coalgebra, \( C \), over an operad, \( \mathscr{V} \), is a sequence
of maps \[
f_{n}:\mathscr{V}\otimes C\to C^{n}\]
 for all \( n>0 \), where \( f_{n} \) is \( RS_{n} \)-equivariant
or maps (via the \emph{adjoint representation}):\[
g_{n}:C\to \homzs{n}(\mathscr{V}_{n},C^{n})\]
 This latter description of coalgebras (via adjoint maps) is frequently
more useful for our purposes than the previous one. In the case where
\( \mathscr{V} \) is \emph{unital}, we write\[
\homzs{0}(\mathscr{V}_{0},C^{0})=\ring \]
and identify the adjoint structure map with the augmentation of \( C \)\[
g_{0}=\varepsilon :C\to \ring =\homzs{0}(\mathscr{V}_{0},C^{0})\]

These adjoint maps are related in the sense that they fit into commutative
diagrams:  \begin{equation}\xymatrix@C-20pt{ {C}\ar[r]^{g_n}\ar[dd]_{g_{n+m-1}} & {\homzs{n}(\mathscr{V}_n,C^n)}\ar[d]^{\homz(1,\underbrace{\scriptstyle1\otimes\cdots\otimes g_m\otimes\cdots\otimes 1}_{i^{\mathrm{th}}\mathrm{~position }})} \\ {} & {\homzs{n}(\mathscr{V}_n,C^{i-1}\otimes\homzs{m}(\mathscr{V}_{m},C^{m})\otimes C^{n-i}) }\ar[d]^{\iota} \\ {\homzs{n+m-1}(\mathscr{V}_{n+m-1},C^{n+m-1})}\ar[r]_-{\homz(\circ_i,1)} & {\homz(\mathscr{V}_n\otimes\mathscr{V}_m,C^{n+m-1}) } }\label{dia:coalgdef}\end{equation}
for all \( m,n>0 \) and all \( 1\le i\le n \), where \( \iota  \)
is the composite \begin{equation}\xymatrix@R-10pt{{\homzs{n}(\mathscr{V}_n,C^{i-1}\otimes\homzs{m}(\mathscr{V}_{m},C^{m})\otimes C^{n-i})} \ar@{=}[d] \\ {\homzs{n}(\mathscr{V}_n,\homz(\ring,C^{i-1})\otimes\homzs{m}(\mathscr{V}_{m},C^{m})\otimes \homz(\ring,C^{n-i}))} \ar@{=}[d] \\ {\homz(\mathscr{V}_n,\homz(\mathscr{V}_{m},C^{m+n-1}))}\ar@{^{(}->}[d] \\ {\homz(\mathscr{V}_n\otimes\mathscr{V}_m,C^{n+m-1}) }}\label{eq:iotadef}\end{equation} 

In other words: The abstract composition-operations in \( \mathscr{V} \)
exactly correspond to compositions of maps in \( \{\homz (C,C^{n})\} \). 

The following is clear:
\end{rem}
\begin{prop}
\label{prop:everycomodule} Every chain complex is trivially a coalgebra
over its own coendomorphism operad. 
\end{prop}

\subsection{Examples}

\begin{example}
Coassociative coalgebras are precisely the coalgebras over \( \mathfrak{S}_{0} \)
(see \ref{def:mathfrakS0}). 
\end{example}
\begin{defn}
\label{def:coassoc}\( \mathbf{Cocommut} \) is an operad defined
to have one basis element \( \{b_{i}\} \) for all integers \( i\ge 0 \).
Here the rank of \( b_{i} \) is \( i \) and the degree is 0 and
the these elements satisfy the composition-law: \( \gamma (b_{n}\otimes b_{k_{1}}\otimes \cdots \otimes b_{k_{n}})=b_{K} \),
where \( K=\sum _{i=1}^{n}k_{i} \). The differential of this operad
is identically zero. The symmetric-group actions are trivial.
\end{defn}
\begin{example}
Coassociative commutative coalgebras are the coalgebras over \( \mathbf{Cocommut} \).
\end{example}
The following example has many topological applications

\begin{example}
Coalgebras over the operad \( \mathfrak{S} \), defined in \ref{def:sfrakfirstmention},
are chain-complexes equipped with a coassociative coproduct and Steenrod
operations for all primes (see \cite{Smith:2000}).
\end{example}

\section{\label{sec:generalconstruction}The general construction}

We begin by defining

\begin{defn}
\label{def:addtree}Let \( n\ge 1 \) be an integer and let \( k \)
be \( 0 \) or \( 1 \). Define \( \mathrf {P}_{k}(n) \) to be the
set of sequences \( \{u_{1},\dots ,u_{m}\} \) of elements each of
which is either a \( \bullet  \)-symbol or an integer \( \ge k \)
and such that\begin{equation}
\label{eq:addtreeequation}
\sum _{j=1}^{m}u_{j}=n
\end{equation}
where \( \bullet =1 \) for the purpose of computing this sum.

Given a sequence \( \mathbf{u}\in \mathrf {P}_{k}(n) \), let \( \slength{\mathbf{u}}=m \),
the length of the sequence.
\end{defn}
\begin{rem}
Note that the set \( \mathrf {P}_{1}(n) \) is finite and for any
\( \mathbf{u}\in \mathrf {P}_{k}(n) \) \( \slength{\mathbf{u}}\le n \).
By contrast, \( \mathrf {P}_{0}(n) \) is always infinite.
\end{rem}
\begin{defn}
\label{def:bigdiag}Let \( \mathscr{V} \) be an operad, let \( C \)
be a \( \ring  \)-free DG module and set\[
k=\left\{ \begin{array}{ll}
0 & \mathrm{if}\, \, \mathscr{V}\, \, \mathrm{is}\, \mathrm{unital}\\
1 & \mathrm{otherwise}
\end{array}\right. \]
Now define \[
KC=C\oplus \prod _{n\ge k}\homzs{n}(\mathscr{V}_{n},C^{n})\]
where \( \homzs{0}(\mathscr{V}_{0},C^{0})=\ring  \) in the unital
case.

Consider the diagram \begin{equation} \xymatrix@C+30pt{{} & {\displaystyle\prod_{m\ge k}\homzs{m}(\mathscr{V}_m,(KC)^m)} \ar@{^{(}->}[d]^-{y}\\ {KC\vphantom{\displaystyle\prod_{\substack{ n\ge k\\ \mathbf{u}\in\mathrf{P}_k(n) }}}}\ar@<1.8ex>[r]_-{0\oplus\prod_{n\ge k}c_n} & {\displaystyle\!\!\!\!\!\prod_{\substack{ n\ge k\\ \mathbf{u}\in\mathrf{P}_k(n) }}\!\!\!\!\!\homz(\mathscr{V}_{\slength{\mathbf{u}}}\otimes \mathscr{V}_{u_{1}}\otimes\cdots\otimes \mathscr{V}_{u_{\slength{\mathbf{u}}}},C^n)}} \label{eq:bigdiag}\end{equation}where
\begin{enumerate}
\item the \( c_{n} \) are defined by \[c_n=\!\!\!\!\!\prod_{ \mathbf{u}\in\mathrf{P}_k(n)}\!\!\!\!\!\homz(\gamma_{\mathbf{u}},1)\]
with\[
\homz (\gamma _{\mathbf{u}},1):\homzs{n}(\mathscr{V}_{n},C^{n})\to \homz (\mathscr{V}_{\slength{\mathbf{u}}}\otimes \mathscr{V}_{u_{1}}\otimes \cdots \otimes \mathscr{V}_{u_{\slength{\mathbf{u}}}},C^{n})\]
 the dual of the generalized structure map \[
\gamma _{\mathbf{u}}:\mathscr{V}_{\slength{\mathbf{u}}}\otimes \mathscr{V}_{u_{1}}\otimes \cdots \otimes \mathscr{V}_{u_{\slength{\mathbf{u}}}}\to \mathscr{V}_{n}\]
 from definition~\ref{def:generalizedcomps}. We assume that \( \mathscr{V}_{\bullet }=\ring  \)
and \( C^{\bullet }=C \) so that \( \homzs{\bullet }(\mathscr{V}_{\bullet },C^{\bullet })=C \).
\item \( y=\prod _{m\ge k}y_{m} \) and the maps \[y_m: \homzs{m}(\mathscr{V}_m,(KC)^m)\to \!\!\!\!\!\prod_{\substack{n\ge k\\ \mathbf{u}\in\mathrf{P}_k(n)}}\homz(\mathscr{V}_{\slength{\mathbf{u}}}\otimes\mathscr{V}_{u_{1}}\otimes\cdots\otimes \mathscr{V}_{u_{\slength{\mathbf{u}}}},C^n)\]
map the factor \[
\homzs{m}(\mathscr{V}_{m},\bigotimes _{j=1}^{m}L_{j})\subset \homzs{m}(\mathscr{V}_{m},(KC)^{m})\]
with \( L_{j}=\homzs{u_{j}}(\mathscr{V}_{u_{j}},C^{u_{j}}) \) via
the map induced by the associativity of the \( \mathrm{Hom} \) and
\( \otimes  \) functors. 
\end{enumerate}
A submodule \( M\subseteq KC \) is called \emph{\( \mathscr{V} \)-closed}
if\[
\left( 0\oplus \prod _{n\ge k}c_{n}\right) (M)\subseteq y\left( \prod _{n\ge k}\homzs{n}(\mathscr{V}_{n},M^{n})\right) \]

\end{defn}
Now we take stock of the terms in diagram~\ref{eq:bigdiag}. 

\begin{lem}
\label{lemma:lvconstruction}Let \( \mathscr{V} \) be an operad and
let \( k=0 \) if \( \mathscr{V} \) is unital and \( 1 \) otherwise.
Under the hypotheses of definition~\ref{def:bigdiag}, if \( C \)
is a DG module over \( \ring  \), set

\begin{eqnarray*}
L_{1}C & = & KC\\
L_{n}C & = & \left( 0\oplus \prod _{n\ge k}c_{n}\right) ^{-1}\prod _{m\ge k}y_{m}\left( \homzs{m}(\mathscr{V}_{m},(L_{n-1}C)^{m})\right) 
\end{eqnarray*}
Then\begin{equation}
\label{eq:lvconstruction}
L_{\mathscr{V}}C=\bigcap _{n=1}^{\infty }L_{n}C
\end{equation}
--- the maximal \( \mathscr{V} \)-closed submodule of \( KC \) (in
the notation of definition~\ref{def:bigdiag}) --- is a coalgebra
over \( \mathscr{V} \) with coproduct given by\begin{equation}
\label{th:cofreecoalgformula}
g=y^{-1}\circ \left( 0\oplus \prod _{n\ge k}c_{n}\right) :L_{\mathscr{V}}C\to \prod _{n\ge k}\homzs{n}(\mathscr{V}_{n},(L_{\mathscr{V}}C)^{n})
\end{equation}

\end{lem}
\begin{rem}
See appendix~\ref{app:proofcoalgebra} for the proof.
\end{rem}
\begin{lem}
\label{lem:classifyingmap}Let \( \mathscr{V} \) be an operad and
let \( k=0 \) if \( \mathscr{V} \) is unital and \( 1 \) otherwise.
Given a coalgebra \( D \) over \( \mathscr{V} \) with adjoint structure
maps\[
d_{n}:D\to \homzs{n}(\mathscr{V}_{n},D^{n})\]
any morphism of DG-modules\[
f:D\to C\]
induces a map\[
\hat{f}=f\oplus \prod _{n\ge k}^{\infty }\homzs{n}(1,f^{n})\circ d_{n}:D\to C\oplus \prod _{n\ge k}\homzs{n}(\mathscr{V}_{n},C^{n})\]
whose image lies in\[
L_{\mathscr{V}}C\subseteq C\oplus \prod _{n\ge k}\homzs{n}(\mathscr{V}_{n},C^{n})\]
as defined in lemma~\ref{lemma:lvconstruction}. Furthermore, \( \hat{f} \)
is a morphism of \( \mathscr{V} \)-coalgebras. 
\end{lem}
\begin{rem}
In the unital case, the augmentation \( L_{\mathscr{V}}C\to \ring  \)
is induced by projection to the factor \( \homzs{0}(\mathscr{V}_{0},C^{0})=\ring  \).
\end{rem}
\begin{proof}
We prove the claim when \( C=D \) and use the functoriality of \( L_{\mathscr{V}}C \)
with respect to morphisms of \( C \) to conclude it in the general
case. In this case \( f=\mathrm{id}:D\to D \) and \( \hat{f}=d \). 

We claim that the diagram \begin{equation}\smaller\xymatrix@C+6pt{ {D}\ar[r]^-{d}\ar[dd]_{\hat{f}}\ar[rrdd]|-{y\circ\prod_{n\ge k}\homz(1,d^n)\circ d} & {\prod_{n\ge k}\homzs{n}(\mathscr{V}_n,D^n)}\ar[r]^-{\prod_{n\ge k}\homz(1,{\hat{f}}^n)} & {\prod_{n\ge k}\homzs{n}(\mathscr{V}_n,(L_{\mathscr{V}}D)^n)}\ar[dd]^{y} \\ {} & {} {} \\ {L_{\mathscr{V}}D\vphantom{\displaystyle\prod_{\substack{ n\ge k \\ \mathbf{u}\in\mathrf{P}_k(n) }}}}\ar@<1.9ex>[rr]_-{0\oplus\prod_{n\ge k}c_n} & {} & {\displaystyle\!\!\!\!\!\prod_{\substack{ n\ge k \\ \mathbf{u}\in\mathrf{P}_k(n) }}\!\!\!\!\!\homz(V_{\slength{\mathbf{u}}}\otimes V_{u_{1}}\otimes\cdots\otimes V_{u_{\slength{\mathbf{u}}}},C^n)} }\end{equation}
commutes, where \( c_{n} \), \( y \), and \( y_{n} \) are as defined
in definition~\ref{def:bigdiag} so that the lower row and right
column are the same as diagram~\ref{eq:bigdiag}. Clearly, the upper
sub-triangle of this diagram commutes since \( \hat{f}=d \). On the
other hand, the lower sub-triangle also clearly commutes by the definition
of \( c_{n} \) and the fact that \( D \) is a \( \mathscr{V} \)-coalgebra.
It follows that the entire diagram commutes. But this implies that
\( \im \hat{f}\subseteq \prod _{n\ge k}\homzs{n}(\mathscr{V}_{n},D^{n})=KD \)
(in the notation of definition~\ref{def:bigdiag}) satisfies the
condition that\[
\left( 0\oplus \prod _{n\ge k}c_{n}\right) (\im \hat{f})\subseteq y\left( \prod _{n\ge k}\homzs{n}(\mathscr{V}_{n},(\im \hat{f})^{n}\right) \]
so that \( \im \hat{f} \) is \( \mathscr{V} \)-closed --- see definition~\ref{def:bigdiag}.
It follows that \( \im \hat{f}\subseteq L_{\mathscr{V}}D\subseteq KD \)
since \( L_{\mathscr{V}}D \) is \emph{maximal} with respect to this
property (see lemma~\ref{lemma:lvconstruction}).

This implies both of the statements of this lemma.
\end{proof}
\begin{thm}
\label{th:cofreegeneral}Let \( D \) a coalgebra over the operad
\( \mathscr{V} \) with adjoint structure maps\[
d_{n}:D\to \homzs{n}(\mathscr{V}_{n},D^{n})\]
and let \[
f:D\to C\]
be any morphism of DG-modules. Then the coalgebra morphism\[
\hat{f}:D\to L_{\mathscr{V}}C\]
 defined in lemma~\ref{lem:classifyingmap} is the unique coalgebra
morphism that makes the diagram \begin{equation}\xymatrix{{D} \ar[r]^-{\hat{f}}\ar[rd]_{f}& {L_{\mathscr{V}}C}\ar[d]^{\epsilon} \\ & {C}}\label{dia:coveringfhat}\end{equation} commute.
Consequently \( L_{\mathscr{V}}C \) is the cofree coalgebra over
\( \mathscr{V} \) cogenerated by \( C \). The cogeneration map (see
definition~\ref{def:cofreecoalgebra}) \( \epsilon :L_{\mathscr{V}}C\to C \)
is projection to the first direct summand.
\end{thm}
\begin{proof}
Let \( k=0 \) if \( \mathscr{V} \) is unital and \( 1 \) otherwise.
It is very easy to see that diagram~\ref{dia:coveringfhat} commutes
with \( \hat{f} \) as defined in lemma~\ref{lem:classifyingmap}.
Suppose that\[
g:D\to L_{\mathscr{V}}C\]
is another coalgebra morphism that makes diagram~\ref{dia:coveringfhat}
commute. We claim that \( g \) must coincide with \( \hat{f} \).
The component \[
\homz (\gamma _{\{\bullet ,\dots ,\bullet \}},1):L_{\mathscr{V}}C\to \prod _{n\ge k}\homzs{n}(\mathscr{V}_{n},(L_{\mathscr{V}}C)^{n})\]
isomorphically maps \[
\homzs{n}(\mathscr{V}_{n},C^{n})\]
to the direct summand\[
\homzs{n}(\mathscr{V}_{n},C^{n})\subset \homzs{n}(\mathscr{V}_{n},(L_{\mathscr{V}}C)^{n})\]
For \( g \) to be a coalgebra morphism, we \emph{must} have (at least)\[
\homz (1,g^{n})\circ d_{n}=\homz (\gamma _{\{\bullet ,\dots ,\bullet \}},1)\circ g\]
for all \( n\ge k \). This requirement, however, \emph{forces} \( g=\hat{f} \).

Lemma~\ref{lem:classifyingmap} and the argument above verify all
of the conditions in definition~\ref{def:cofreecoalgebra}.
\end{proof}

\section{Special cases\label{sec:special}}

\subsection{Dimension restrictions}

Now we address the issue of our cofree coalgebra extending into negative
dimensions. We need the following definition first:

\begin{defn}
\label{def:truncation}If \( E \) is a chain-complex, and \( t \)
is an integer, let \( E^{\dimlimiter t} \) denote the chain-complex
defined by\[
E^{\dimlimiter t}_{i}=\left\{ \begin{array}{ll}
0 & \mathrm{if}\, i\le t\\
\ker \partial _{t+1}:E_{t+1}\to E_{t} & \mathrm{if}\, i=t+1\\
E_{i} & \mathrm{if}\, i>t+1
\end{array}\right. \]

\end{defn}
\begin{cor}
\label{th:cofreecoalgebranonnegative}If \( C \) is a chain-complex
concentrated in nonnegative dimensions and \( \mathscr{V} \) is an
operad, then there exists a sub-\( \mathscr{V} \)-coalgebra \[
M_{\mathscr{V}}C\subset L_{\mathscr{V}}C\]
such that
\begin{enumerate}
\item as a chain-complex, \( M_{\mathscr{V}}C \) is concentrated in nonnegative
dimensions,
\item for any \( \mathscr{V} \)-coalgebra, \( D \), concentrated in nonnegative
dimensions, the image of the classifying map \[
\hat{f}:D\to L_{\mathscr{V}}C\]
lies in \( M_{\mathscr{V}}C\subset L_{\mathscr{V}}C \).
\end{enumerate}
In addition, \( M_{\mathscr{V}}C \) can be defined inductively as
follows: Let \( M_{0}=(L_{\mathscr{V}}C)^{\dimlimiter -1} \) (see
\ref{def:truncation}) with structure map\[
\delta _{0}:M_{0}\to \prod _{n>0}\homzs{n}(\mathscr{V}_{n},(L_{\mathscr{V}}C)^{n})=Q_{-1}\]
Now define\[
M_{i+1}=\delta _{i}^{-1}\left( \delta _{i}(M_{i})\cap Q_{i}\right) \subseteq \delta _{i}^{-1}Q_{i-1}\]
 with structure map\[
\delta _{i+1}=\delta _{i}|M_{i+1}:M_{i+1}\to Q_{i}\]
where\[
Q_{i}=\prod _{n>0}\homzs{n}(\mathscr{V}_{n},M_{i}^{n})^{\dimlimiter -1}\]
Then \[
M_{\mathscr{V}}C=\bigcap _{i=0}^{\infty }M_{i}\]

\end{cor}
\begin{rem}
Our definition of \( M_{\mathscr{V}}C \) is simply that of the maximal
sub-coalgebra of \( L_{\mathscr{V}}C \) contained within \( L_{\mathscr{V}}C^{\dimlimiter -1} \)
--- i.e., the maximal sub-coalgebra in \emph{nonnegative} dimensions.
\end{rem}

\subsection{The pointed irreducible case}

We define the pointed irreducible coalgebras over an operad in a way
that extends the conventional definition in \cite{Sweedler:1969}:

\begin{defn}
\label{def:pointedirreducible} Given a coalgebra over a unital operad
\( \mathscr{V} \) with adjoint structure map\[
a_{n}:C\to \homzs{n}(\mathscr{V}_{n},C^{n})\]
an \emph{element} \( c\in C \) is called \emph{group-like} if \( a_{n}(c)=f_{n}(c^{n}) \)
for all \( n>0 \). Here \( c^{n}\in C^{n} \) is the \( n \)-fold
\( \ring  \)-tensor product, \[
f_{n}=\homz (\epsilon _{n},1):\homz (\ring ,C^{n})=C^{n}\to \homzs{n}(\mathscr{V}_{n},C^{n})\]
and \( \epsilon _{n}:\mathscr{V}_{n}\to \ring  \) is the augmentation
(which exists by \ref{def:unitaloperad}). 

A coalgebra \( C \) over a unital operad \( \mathscr{V} \) is called
\emph{pointed} if it has a \emph{unique} group-like element (denoted
\( 1 \)), and \emph{pointed irreducible} if the intersection of any
two sub-coalgebras contains this unique group-like element.
\end{defn}
\begin{rem}
Note that a group-like element generates a sub \( \mathscr{V} \)-coalgebra
of \( C \) and must lie in dimension \( 0 \).

Although this definition seems contrived, it arises in {}``nature'':
The chain-complex of a pointed, simply-simply connected reduced simplicial
set is pointed irreducible over the operad \( \mathfrak{S} \). In
this case, the operad action encodes the effect on the chain level
of all Steenrod operations.
\end{rem}
Note that our cofree coalgebra in theorem~\ref{th:cofreecoalgformula}
is pointed since it has the sub-coalgebra \( \ring  \). It is \emph{not}
irreducible since the \emph{null} submodule, \( C \) (on which the
coproduct vanishes identically), is a sub-coalgebra whose intersection
with \( \ring  \) is \( 0 \). We conclude that:

\begin{lem}
\label{lem:irrednodead}Let \( C \) be a pointed irreducible coalgebra
over a unital operad \( \mathscr{V} \). Then the adjoint structure
map\[
C\to \prod _{n\ge 0}\homzs{n}(\mathscr{V}_{n},C^{n})\]
 is injective.
\end{lem}
The existence of units of operads, and the associativity relations
imply that

\begin{lem}
\label{lem:a1mapsplit}Let \( C \) be a coalgebra over an operad
\( \mathscr{V} \) with the property that the adjoint structure map\[
\prod _{n\ge 1}a_{n}:C\to \prod _{n\ge 1}\homzs{n}(\mathscr{V}_{n},C^{n})\]
 is injective. Then the adjoint structure map\[
a_{1}:C\to \homz (\mathscr{V}_{1},C)\]
 is naturally split by\[
\homz (\eta _{1},1):\homz (\mathscr{V}_{1},C)\to \homz (\ring ,C)=C\]
where \( \eta _{1}:\ring \to \mathscr{V}_{1} \) is the unit.
\end{lem}
\begin{rem}
In general, the unit \( \eta _{1}\in \mathscr{V} \) maps under the
structure map\[
s:\mathscr{V}\to \coend (C)\]
to a unit of \( \im s \) --- a \emph{sub}-operad of \( \coend (C) \).
We show that \( s(\eta _{1}) \) is \( 1:C\to C\in \coend (C)_{1} \).
\end{rem}
\begin{proof}
Consider the endomorphism\[
e=\homz (\eta _{1},C)\circ a_{1}:C\to C\]
The operad identities imply that the diagram \[\xymatrix@C+20pt{{C}\ar[r]^-{\prod_{n\ge 1}a_n}\ar[d]_{e} & {\prod_{n\ge 1}\homzs{n}(\mathscr{V}_n,C^n)} \\ {C}\ar[ru]_-{\prod_{n\ge 1}a_n} & {}}\]
commutes since \( \eta _{1} \) is a unit of the operad and \( \homz (\eta _{1},C)\circ a_{1} \)
must preserve the coproduct structure (acting, effectively, as the
\emph{identity map}). 

It follows that \( e^{2}=e \) and that \( \ker e\subseteq \ker \prod _{n\ge 1}a_{n} \).
The hypotheses imply that \( \ker e=0 \) and we claim that \( e^{2}=e\Rightarrow \im e=C \).
Otherwise, suppose that \( x\in C\setminus \im e \). Then \( e(x-e(x))=0 \)
so \( x-e(x)\in \ker e \), which is a contradiction. The conclusion
follows. 
\end{proof}
\begin{prop}
Let \( D \) be a pointed, irreducible coalgebra over a unital operad
\( \mathscr{V} \). Then the augmentation map\[
\varepsilon :D\to \ring \]
is naturally split and any morphism of pointed, irreducible coalgebras
\[
f:D_{1}\to D_{2}\]
 is of the form\[
1\oplus \bar{f}:D_{1}=\ring \oplus \ker \varepsilon _{D_{1}}\to D_{2}=\ring \oplus \ker \varepsilon _{D_{2}}\]
where \( \varepsilon _{i}:D_{i}\to \ring  \), \( i=1,2 \) are the
augmentations.
\end{prop}
\begin{proof}
The definition (\ref{def:pointedirreducible}) of the sub-coalgebra
\( \ring \subseteq D_{i} \) is stated in an invariant way, so that
any coalgebra morphism must preserve it.
\end{proof}
Our result is:

\begin{thm}
\label{th:pointedirreduciblecofree}If \( C \) is a chain-complex
and \( \mathscr{V} \) is a unital operad, define\[
P_{\mathscr{V}}C=L_{\mathscr{V}}C/C\]
 (see theorem~\ref{th:cofreecoalgformula}) with the induced quotient
structure map.

Then \( P_{\mathscr{V}}C \) is a pointed, irreducible coalgebra over
\( \mathscr{V} \). Given any pointed, irreducible coalgebra \( D \)
over \( \mathscr{V} \) with adjoint structure maps\[
d_{n}:D\to \homzs{n}(\mathscr{V}_{n},D^{n})\]
and augmentation\[
\varepsilon :D\to \ring \]
any morphism of DG-modules\[
f:\ker \varepsilon \to C\]
extends to a unique morphism of pointed, irreducible coalgebras over
\( \mathscr{V} \)\[
1\oplus \hat{f}:\ring \oplus \ker \varepsilon \to P_{\mathscr{V}}C\]
where\[
\hat{f}=1\oplus \prod _{n=1}^{\infty }\homzs{n}(1,f^{n})\circ d_{n}:D\to P_{\mathscr{V}}C\]

If \( p_{C}:P_{\mathscr{V}}C\to \homz (\mathscr{V}_{1},C) \) is projection
to the first factor, and \( \homz (\eta _{1},1):\homz (\mathscr{V}_{1},C)\to C \)
is the splitting map defined in \ref{lem:a1mapsplit}, then the diagram
\begin{equation}\xymatrix{{D} \ar[r]^-{\hat{f}}\ar[rd]_{f}& {P_{\mathscr{V}}C}\ar[d]^{\homz(\eta_1,1)\circ p_C} \\ & {C}}\label{dia:cogeneration}\end{equation}
commutes.

In particular, \( P_{\mathscr{V}}C \) is the cofree pointed irreducible
coalgebra over \( \mathscr{V} \) with cogeneration map \( \homz (\eta ,1)\circ p_{C} \)
(see definition~\ref{def:cofreecoalgebra}).
\end{thm}
\begin{rem}
Roughly speaking, \( P_{\mathscr{V}}C \) is an analogue to the \emph{Shuffle
Coalgebra} defined in \cite[chapter 11]{Sweedler:1969}. With one
extra condition on the operad \( \mathscr{V} \), this becomes a generalization
of the Shuffle Coalgebra.
\end{rem}
\begin{proof}
Since the kernel of the the structure map of \( D \) vanishes\[
\im \hat{f}\cap C=0\]
 so that \( \im \hat{f} \) is mapped isomorphically by the projection
\( L_{\mathscr{V}}C\to P_{\mathscr{V}}C \). 

It is first necessary to show that \( \homz (\eta _{1},1)\circ p_{C}:\homz (\mathscr{V}_{1},C)\to C \)
can serve as the cogeneration map, i.e., that diagram~\ref{dia:cogeneration}
commutes.

This conclusion follows from the commutativity of the diagram \[\xymatrix@+20pt{ {D}\ar[r]^{d}\ar@/_/[rd]_{=}\ar@/^2pc/[rr]^{\hat{f}} & {P_{\mathscr{V}}D}\ar[r]^{P_{\mathscr{V}}f}\ar[d]|-{\homz(\eta_1,1_D)\circ p_D} & {P_{\mathscr{V}}C}\ar[d]^{\homz(\eta_1,1_C)\circ p_C} \\ {} & {D}\ar[r]_{f} & {C}}\] where
\( d:D\to P_{\mathscr{V}}C \) is the canonical classifying map of
\( D \).

The upper (curved) triangle commutes by the definition of \( \hat{f} \),
the lower left triangle by the fact that \( \homz (\eta _{1},1) \)
splits the classifying map. The lower right square commutes by functoriality
of \( P_{\mathscr{V}}\ast  \). 

We must also show that \( P_{\mathscr{V}}C \) is pointed irreducible.
The sub-coalgebra generated by \( 1\in \ring =\homzs{0}(\mathscr{V}_{0},C^{0}) \)
is group-like. 

\emph{Claim: If \( x\in P_{\mathscr{V}}C \) is an arbitrary element,
its coproduct in \( \homzs{N}(\mathscr{V},P_{\mathscr{V}}C^{N}) \)
for \( N \) sufficiently large, contains factors of \( 1\in \ring \subset P_{\mathscr{V}}C \). }

This follows from the fact that \( \mathbf{u}\in \mathrf {P}_{0}(n) \)
\emph{must} have terms \( u_{i}=0 \) for \( N=\slength{\mathbf{u}}>n \)
--- see \ref{def:addtree} with \( k=0 \). 

It follows that \emph{every} sub-coalgebra of \( P_{\mathscr{V}}C \)
must contain \( 1 \) so that \( \ring  \) is the \emph{unique} sub-coalgebra
of \( P_{\mathscr{V}}C \) generated by a group-like element. This
implies that \( P_{\mathscr{V}}C \) is pointed irreducible.

The statement about any pointed irreducible coalgebra mapping to \( P_{\mathscr{V}}C \)
follows from lemma~\ref{lem:classifyingmap}.
\end{proof}
\begin{defn}
\label{def:kconnectedirred}Let \( C \) be a pointed irreducible
\( \mathscr{V} \)-coalgebra with augmentation\[
\varepsilon :C\to \ring \]
If \( t \) is some integer, we say that \( C \) is \( t \)-\emph{reduced}
if \[
(\ker \varepsilon )_{i}=0\]
for all \( i\le t \).
\end{defn}
\begin{rem}
If \( t\ge 1 \), the chain complex of a \( t \)-reduced simplicial
set (see \cite[p.~170]{Goerss-Jardine}) is naturally a \( t \)-reduced
pointed, irreducible coalgebra over \( \mathfrak{S} \). The case
where \( t\le 0 \) also occurs in topology in the study of spectra.
\end{rem}
We conclude this section with a variation of \ref{th:cofreecoalgebranonnegative}. 

\begin{prop}
\label{prop:pointedirreducnonnegative}If \( t \) is an integer and
\( C \) is a chain-complex concentrated in dimensions \( >t \),
and \( \mathscr{V} \) is a unital operad, let \( P_{\mathscr{V}}C \)
be the pointed, irreducible coalgebra over \( \mathscr{V} \) defined
in \ref{th:pointedirreduciblecofree}. There exists a sub-coalgebra,
\[
\mathrf {F}^{\dimlimiter t}_{\mathscr{V}}C\subset P_{\mathscr{V}}C\]
such that
\begin{enumerate}
\item \( \mathrf {F}^{\dimlimiter t}_{\mathscr{V}}C \) is a \( t \)-reduced
pointed irreducible coalgebra over \( \mathscr{V} \),
\item for any pointed, irreducible \( t \)-reduced \( \mathscr{V} \)-coalgebra,
\( D \), the image of the classifying map \[
1\oplus \hat{f}:D\to P_{\mathscr{V}}C\]
lies in \( \mathrf {F}^{\dimlimiter t}_{\mathscr{V}}C\subset P_{\mathscr{V}}C \).
\end{enumerate}
In addition, \( \mathrf {F}^{\dimlimiter t}_{\mathscr{V}}C \) can
be defined inductively as follows: Let \( Y_{0}=\ring \oplus (P_{\mathscr{V}}C)^{\dimlimiter t} \)
(see \ref{def:truncation} for the definition of \( (\ast )^{\dimlimiter t} \))
with structure map\[
\alpha _{0}:Y_{0}\to \ring \oplus \prod _{n>0}\homzs{n}(\mathscr{V}_{n},(P_{\mathscr{V}}C)^{n})=Z_{-1}\]
Now define \begin{equation}
\label{eq:intersect}
Y_{i+1}=\alpha _{i}^{-1}\left( \alpha _{i}(Y_{i})\cap Z_{i}\right) \subseteq \alpha _{i}^{-1}Z_{i-1}
\end{equation}
with structure map\[
\alpha _{i+1}=\alpha _{i}|Y_{i+1}:Y_{i+1}\to Z_{i}\]
where\[
Z_{i}=\ring \oplus \prod _{n>0}\homzs{n}(\mathscr{V}_{n},Y_{i}^{n})^{\dimlimiter t}\]
Then \[
\mathrf {F}^{\dimlimiter t}_{\mathscr{V}}C=\bigcap _{i=0}^{\infty }Y_{i}\]

\end{prop}
\begin{rem}
Our definition of \( \mathrf {F}^{\dimlimiter t}_{\mathscr{V}}C \)
is simply that of the maximal sub-coalgebra of \( P_{\mathscr{V}}C \)
contained within \( \ring \oplus P_{\mathscr{V}}C^{\dimlimiter t} \).
\end{rem}
\begin{example}
For example, let \( \mathscr{V}=\mathfrak{S}_{0} \) --- the operad
whose coalgebras are coassociative coalgebras. Let \( C \) be a chain-complex
concentrated in positive dimensions. Since the operad is concentrated
in dimension \( 0 \) the {}``natural'' coproduct given in \ref{th:pointedirreduciblecofree}
does \emph{not} go into negative dimensions when applied to \( \ring \oplus \prod _{n>0}\homzs{{n}}(\mathscr{V}_{n},C^{n})^{\dimlimiter 0} \)
so \( M_{n}C=\homzs{n}(\mathscr{V}_{n},C^{n})^{\dimlimiter 0}=\homzs{n}(\mathscr{V}_{n},C^{n}) \)
for all \( n>0 \) and \begin{eqnarray*}
\mathrf {F}^{\dimlimiter 0}_{\mathscr{V}}C & = & \ring \oplus \prod _{n>0}\homzs{{n}}(\mathscr{V}_{n},C^{n})^{\dimlimiter 0}\\
 & = & \ring \oplus \bigoplus _{n>0}\homzs{{n}}(\zs{{n}},C^{n})\\
 & = & T(C)
\end{eqnarray*}
the \emph{tensor-algebra} --- the well-known pointed, irreducible
cofree coalgebra used in the \emph{bar construction.} 

The fact that the direct product is of \emph{graded} modules and dimension
considerations imply that, in each dimension, it only has a \emph{finite}
number of nonzero factors. So, in this case, the direct product becomes
a direct sum.
\end{example}
\appendix

\section{Proof of lemma~\ref{lemma:lvconstruction}\label{app:proofcoalgebra}}

As always, \( k=0 \) if \( \mathscr{V} \) is unital and \( 1 \)
otherwise. Note that the coproduct formula, equation~\ref{th:cofreecoalgformula}
is well-defined because the map\[
y=\prod _{m\ge k}y_{m}\]
is injective and\[
(0\oplus \prod _{n\ge k}c_{n})(L_{\mathscr{V}}C)\subseteq y\left( \prod _{n\ge k}\homzs{n}(\mathscr{V}_{n},(L_{\mathscr{V}}C)^{n})\right) \]
 by our construction of \( L_{\mathscr{V}}C \) in equation~\ref{eq:lvconstruction}. 

The basic idea behind this proof is that we \emph{dualize} the argument
used in verifying the defining properties of a free algebra over an
operad in \cite{Kriz-May}. This is complicated by the fact that \( L_{\mathscr{V}}C \)
is not really the dual of a free algebra. The closest thing we have
to this dual is \( KC \) in definition~\ref{def:bigdiag}. But \( L_{\mathscr{V}}C \)
is \emph{contained} in \( KC \), not equal to it. We cannot dualize
the proof that a free \( \mathscr{V} \)-algebra is free, but can
carry out a similar argument with respect to a kind of {}``Hilbert
basis'' of \( L_{\mathscr{V}}C \).

Consider a factor\[
\homzs{\alpha }(\mathscr{V}_{\alpha },C^{\alpha })\subset C\oplus \prod _{n\ge k}\homzs{n}(\mathscr{V}_{n},C^{n})\]
In general\[
\homzs{\alpha }(\mathscr{V}_{\alpha },C^{\alpha })\not \subset L_{\mathscr{V}}C\subset C\oplus \prod _{n\ge k}\homzs{n}(\mathscr{V}_{n},C^{n})\]
but we still have a projection\[
p_{\alpha }:L_{\mathscr{V}}C\to \homzs{\alpha }(\mathscr{V}_{\alpha },C^{\alpha })\]
Let its image be \( K_{\alpha }\subseteq \homzs{\alpha }(\mathscr{V}_{\alpha },C^{\alpha }) \).
We will show that all faces of the diagram in \begin{sidewaysfigure}$$\xymatrix@R+50pt@C-55pt{&{K_\alpha} \ar[rr]^{y_n^{-1}\circ\homz(\gamma_{\mathbf{u}},1)}\ar'[d]|-{y_{n+m-1}^{-1}\circ\homz(\gamma_{\mathbf{w}},1)}[dd] & & {\homzs{n}(\mathscr{V}_n,\bigotimes_{j=1}^nK_{u_j})} \ar[dd]|-(.4){\iota_1\circ\homz(1,1^{i-1}\otimes(y_m^{-1}\homz(\gamma_{\mathbf{V}},1))\otimes1^{n-i})} \\ {L_{\mathscr{V}}C} \ar[ur]^{p_\alpha} \ar[rr]^(.3){g_n} \ar[dd]_{g_{n+m-1}} & & {\homzs{n}\left(\mathscr{V}_n,(L_{\mathscr{V}}C)^n\right)}\ar[ur]|-{\homzs{n}(1,\bigotimes_{j=1}^np_{u_j})} \ar[dd]|-(.3){\iota\circ\homz(1,\underbrace{\scriptstyle1\otimes\cdots\otimes g_m\otimes\cdots\otimes 1}_{i^{\mathrm{th}}\mathrm{~position }})} & \\ &{\homzs{n+m-1}(\mathscr{V}_{n+m-1},\bigotimes_{j=1}^{n+m-1}K_{w_j})} \ar'[r][rr]_-{\homz(\circ_i,1)} & & {\homz(\mathscr{V}_{n}\otimes\mathscr{V}_{m},\bigotimes_{j=1}^{n+m-1}K_{w_j})} \\ {\homzs{n+m-1}\left(\mathscr{V}_{n+m-1},(L_{\mathscr{V}}C)^{n+m-1}\right)} \ar[rr]_-{\homz(\circ_i,1)} \ar[ur]|-{\homz(1,\bigotimes_{j=1}^{n+m-1}p_{w_j})} & &{\homz\left(\mathscr{V}_n\otimes\mathscr{V}_m,(L_{\mathscr{V}}C)^{n+m-1}\right)} \ar[ur]_{\quad\homz(1\otimes1,\bigotimes_{j=1}^{n+m-1}p_{w_j})}& }$$\caption{}\label{fig:bigcubediagram}\end{sidewaysfigure} figure~\ref{fig:bigcubediagram}
other than the front face commute for all \( \alpha ,n,m \) and \( \mathbf{u}\in \mathrf {P}_{k}(\alpha ) \),
with \( u_{i}\ne \bullet  \) and \( \slength{\mathbf{u}}=n \), \( \mathbf{v}\in \mathrf {P}_{k}(u_{i}) \)
with \( \slength{\mathbf{v}}=m \) and \( \mathbf{w}\in \mathrf {P}_{k}(\alpha ) \)
where \( \mathbf{w} \) is the result of replacing the \( i^{\mathrm{th}} \)
entry of \( \mathbf{u} \) by \( \mathbf{v} \), so \( \slength{\mathbf{w}}=n+m-1 \).
We assume that \( u_{i}\ne \bullet  \) since the coproduct on the
copy of \( C \) represented by \( u_{i}=\bullet  \) would vanish.
Here, \( \iota  \) is the composite \begin{equation}\smaller\xymatrix{ {\homzs{n}(\mathscr{V}_n,L_{\mathscr{V}}C^{i-1}\otimes \homzs{m}(\mathscr{V}_m,L_{\mathscr{V}}C^{m})\otimes L_{\mathscr{V}}C^{n-i})} \ar@{=}[d] \\ {\homzs{n}(\mathscr{V}_n,\homz(\ring,L_{\mathscr{V}}C^{i-1})\otimes \homzs{m}(\mathscr{V}_m,L_{\mathscr{V}}C^{m})\otimes \homz(\ring,L_{\mathscr{V}}C^{n-i}))} \ar[d] \\ {\homz(\mathscr{V}_n,\homz(\mathscr{V}_m,L_{\mathscr{V}}C^{n+m-1}))} \ar[d] \\ {\homz(\mathscr{V}_n\otimes\mathscr{V}_m,L_{\mathscr{V}}C^{n+m-1})} }\label{dia:iotadef}\end{equation} and
\( \iota _{1} \) is the composite \begin{equation}\xymatrix@R-10pt{ {\homzs{n}(\mathscr{V}_n,\bigotimes_{j=1}^{i-1}K_{u_j}\otimes\homzs{m}(\mathscr{V}_m,\bigotimes_{j=1}^mK_{v_j})\otimes\bigotimes_{j=i+1}^nK_{u_j})} \ar@{=}[d] \\ {\homzs{n}(\mathscr{V}_n,A\otimes\homzs{m}(\mathscr{V}_m,\bigotimes_{j=1}^mK_{v_j})\otimes B)} \ar[d] \\ {\homz(\mathscr{V}_n,\homz(\mathscr{V}_m,\bigotimes_{j=1}^{n+m-1}K_{w_j}))} \ar[d] \\ {\homz(\mathscr{V}_n\otimes\mathscr{V}_m,\bigotimes_{j=1}^{n+m-1}K_{w_j})} }\label{dia:iota1def}\end{equation}
with \( A=\bigotimes _{j=1}^{i-1}\homz (\ring ,K_{u_{k}}) \), \( B=\bigotimes _{j=i+1}^{n}\homz (\ring ,K_{u_{k}}) \).
The top face commutes by the definition of the coproduct of \( L_{\mathscr{V}}C \)
and the fact that the image of the coproduct of an element \( x\in L_{\mathscr{V}}C \)
under \( \homz (1,\bigotimes _{k=1}^{n}p_{u_{k}}) \) \emph{only}
depends on \( p_{\alpha }(x) \) --- since \( \sum _{k=1}^{n}u_{k}=\alpha  \).
This also implies that the left face commutes since the left face
is the same as the top face (for \( g_{n+m-1} \) rather than \( g_{n} \)). 

To see that the right face commutes, note that \( \iota  \) and \( \iota _{1} \)
are very similar --- each term of diagram~\ref{dia:iota1def} projections
to the corresponding term of diagram~\ref{dia:iotadef}. The naturality
of the projection maps and the fact that the top face of the diagram
in figure~\ref{fig:bigcubediagram} commutes implies that the right
face commutes.

Note that definition~\ref{def:addtree} implies that\begin{eqnarray*}
\alpha  & = & \sum _{j=1}^{n}u_{j}\\
 & = & \sum _{j=1}^{n+m-1}w_{j}\\
u_{i} & = & \sum _{j=1}^{m}v_{j}
\end{eqnarray*}

Since elements of \( K_{\alpha } \) are determined by their projections,
the commutativity of all faces of the diagram in figure~\ref{fig:bigcubediagram}
except the front also implies that the \emph{front} face commutes.
This will prove lemma~\ref{lemma:lvconstruction} since it implies
that diagram~\ref{dia:coalgdef} of definition~\ref{def:operadcomodule}
commutes.

The bottom face of the diagram in figure~\ref{fig:bigcubediagram}
commutes by the functoriality of \( \homz (\ast ,\ast ) \).

It remains to prove that the \emph{back} face commutes. To establish
this, we consider the diagram \begin{sidewaysfigure}\smaller$$\xymatrix@R+50pt@C-50pt{&{\homzs{\alpha}(\mathscr{V}_{\alpha},C^{\alpha})} \ar[rr]^{\homz(\gamma_{\mathbf{u}},1)}\ar'[d]|-{\homz(\gamma_{\mathbf{w}},1)}[dd] & & {\homz(\mathscr{V}_n\otimes\bigotimes_{j=1}^n\mathscr{V}_{u_j},C^{\alpha})} \ar[dd]|-(.4){\homz(s\circ1\otimes1^{i-1}\otimes\gamma_{\mathbf{v}}\otimes1^{n-i},1)} \\ {K_\alpha} \ar@{^{(}->}[ur]^{} \ar[rr]^(.3){y_n^{-1}\circ\homz(\gamma_{\mathbf{u}},1)} \ar[dd]_{y_{n+m-1}^{-1}\circ\homz(\gamma_{\mathbf{w}},1)} & & {\homzs{n}(\mathscr{V}_n,\bigotimes_{j=1}^nK_{u_j})}\ar@{^{(}->}[ur]^{y_n} \ar[dd]|-(.3){\iota_1\circ\homz(1,1^{i-1}\otimes(y_m^{-1}\homz(\gamma_{\mathbf{v}},1))\otimes1^{n-i})} & \\ &{\homz(\mathscr{V}_{n+m-1}\otimes\bigotimes_{j=1}^{n+m-1}\mathscr{V}_{w_j},C^{\alpha})} \ar'[r][rr]_-{\homz(\circ_i\otimes1,1)} & & {\homz(\mathscr{V}_{n}\otimes\mathscr{V}_{m}\otimes\bigotimes_{j=1}^{n+m-1}\mathscr{V}_{w_j},C^{\alpha})} \\ {\homzs{n+m-1}(\mathscr{V}_{n+m-1},\bigotimes_{j=1}^{n+m-1}K_{w_j})} \ar[rr]_-{\homz(\circ_i,1)} \ar@{^{(}->}[ur]|-{y_{n+m-1}} & &{\homz(\mathscr{V}_{n}\otimes\mathscr{V}_{m},\bigotimes_{j=1}^{n+m-1}K_{w_j})} \ar@{^{(}->}[ur]_{\iota_2}& }$$\caption{}\label{fig:bigcubediagram2}\end{sidewaysfigure}
in figure~\ref{fig:bigcubediagram2}, where

\begin{eqnarray*}
s:\mathscr{V}_{n}\otimes \mathscr{V}_{m}\otimes \bigotimes _{j=1}^{n+m-1}\mathscr{V}_{w_{j}}= &  & \\
\mathscr{V}_{n}\otimes \mathscr{V}_{m}\otimes \bigotimes _{j=1}^{i-1}\mathscr{V}_{u_{j}}\otimes \left( \bigotimes _{\ell =1}^{u_{i}}\mathscr{V}_{v_{\ell }}\right) \otimes \bigotimes _{j=i+1}^{n}\mathscr{V}_{u_{j}} & \to  & \\
\mathscr{V}_{n}\otimes \bigotimes _{j=1}^{i-1}\mathscr{V}_{u_{j}}\otimes \left( \mathscr{V}_{m}\otimes \bigotimes _{\ell =1}^{u_{i}}\mathscr{V}_{v_{\ell }}\right) \otimes \bigotimes _{j=i+1}^{n}\mathscr{V}_{u_{j}} &  & 
\end{eqnarray*}
is the shuffle map and \( \iota _{2} \) is the composite \[\xymatrix{{\homz(\mathscr{V}_n\otimes\mathscr{V}_m,\bigotimes_{j=1}^{n+m-1}K_{w_j})}\ar@{^{(}->}[d]^{\homz(1,\bigotimes_{j=1}^{n+m-1}y_{w_j})} \\ {\homz(\mathscr{V}_n\otimes\mathscr{V}_m,\bigotimes_{j=1}^{n+m-1}\homzs{w_j}(\mathscr{V}_{w_j},C^{w_j}))}\ar[d] \\ {\homz(\mathscr{V}_n\otimes\mathscr{V}_m,\homz(\bigotimes_{j=1}^{n+m-1}\mathscr{V}_{w_j},C^{\alpha}))}\ar[d] \\ {\homz(\mathscr{V}_n\otimes\mathscr{V}_m\otimes\bigotimes_{j=1}^{n+m-1}\mathscr{V}_{w_j},C^{\alpha})}}\] where
the maps in the lower two rows are the natural associativity maps
for the \( \homz  \)- functor and \( \otimes  \).

Clearly, the left and top faces of the diagram in figure~\ref{fig:bigcubediagram2}
commute. The \emph{bottom} face also commutes because

\begin{enumerate}
\item the maps \( \homz (\circ _{i},1) \) and \( \homz (\circ _{i}\otimes 1,1) \)
only affect the first argument in the \( \homz (\ast ,\ast ) \)-functor
and the other maps in the bottom face only affect the second (so there
is no interactions between them)
\item the remaining maps in that face are composites of natural multilinear
associativity maps like those listed in equation~\ref{eq:specmorphism1}
through \ref{eq:specmorphism4}, so they commute by theorem~\ref{th:linearrigidity}.
\end{enumerate}
The \emph{rear} face commutes because the diagram \[\xymatrix{{\mathscr{V}_{\alpha}} & {\mathscr{V}_n\otimes\bigotimes_{k=1}^n\mathscr{V}_{u_k}}\ar[l]_{\gamma_{\mathbf{u}}} \\ {\mathscr{V}_{n+m-1}\otimes\bigotimes_{k=1}^{n+m-1}\mathscr{V}_{w_k}}\ar[u]^{\gamma_{\mathbf{w}}} & {\mathscr{V}_{n}\otimes\mathscr{V}_{m}\otimes\bigotimes_{k=1}^{n+m-1}\mathscr{V}_{w_k}}\ar[l]^{\circ_i\otimes1}\ar[u]_{1^{i-1}\otimes \gamma_{\mathbf{v}}\otimes1^{n-i}\circ s} }\]
commutes due to the associativity relations for an operad --- see
lemma~\ref{lem:generalizedcompassoc}.

It remains to prove that the \emph{right} face of the diagram in figure~\ref{fig:bigcubediagram2}
commutes. We note that all of the morphisms involved in the right
face are of the type listed in equation~\ref{eq:specmorphism1} through
\ref{eq:specmorphism4} except for \( \gamma _{\mathbf{v}} \) and
invoke theorem~\ref{th:rigidity2}

\section{Multilinear functors\label{app:multilinear}}

In this appendix, we consider multilinear functors on the category
of free \( \ring  \)-modules and show that certain natural transformations
of them must be canonically equal. 

\begin{defn}
\label{def:expressiontree}An \emph{expression tree} is a rooted,
ordered tree whose nodes are labeled with symbols \( \homa  \) and
\( \otimes  \) such that
\begin{enumerate}
\item every node labeled with \( \homa  \) has precisely two children,
\item every node labeled with \( \otimes  \) can have an arbitrary (finite)
number of children,
\item leaf nodes are labeled with \emph{distinct} \( \ring  \)-modules.
\end{enumerate}
Nodes are assigned a quality called \emph{variance} (\emph{covariance}
or \emph{contravariance}) as follows:
\begin{enumerate}
\item The root is covariant.
\item All children of a \( \otimes  \)-node and the right child of a \( \homa  \)-node
have the \emph{same} variance as it.
\item The left child of a \( \homa  \)-node is given the \emph{opposite}
variance.
\end{enumerate}
Two expression-trees are regarded as the \emph{same} if there exists
an isomorphism of ordered trees between them that preserves node-labels.
\end{defn}
\begin{rem}
\label{rem:expressionstree1}For instance,

\begingroup\psfrag{A1}{$A_1$}

\psfrag{A2}{$A_2$}

\psfrag{A3}{$A_3$}

\psfrag{A4}{$A_4$}

\psfrag{A5}{$A_5$}

\psfrag{A6}{$A_6$}

\centerline{\resizebox*{3in}{3in}{\includegraphics{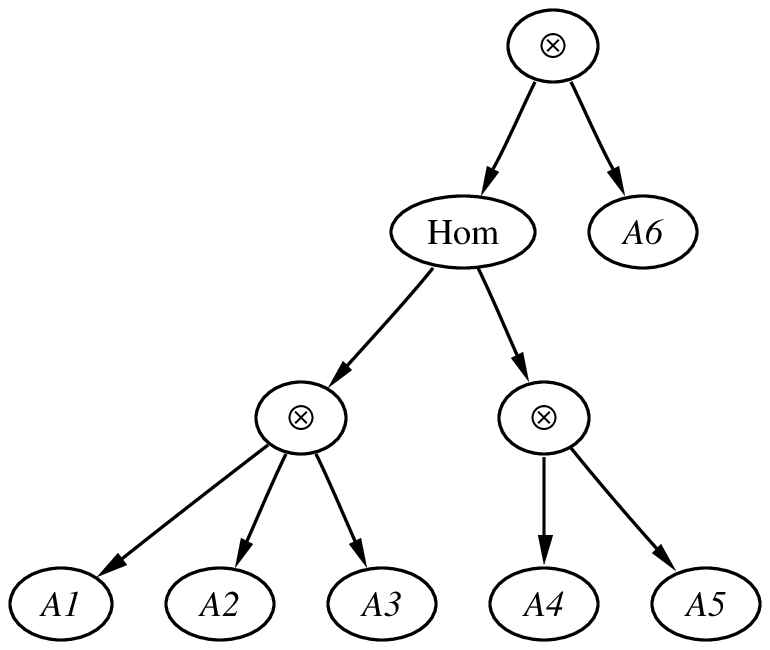}} }\endgroup

\noindent is an expression tree. That expression-trees are \emph{rooted
and ordered} means that:
\begin{enumerate}
\item there is a distinguished node called the \emph{root} that is preserved
by isomorphisms
\item the children of every interior node have a well-defined \emph{ordering}
that is preserved by any isomorphism
\end{enumerate}
\end{rem}
\begin{defn}
\label{def:exptreemodule}Given an expression tree \( T \), let \( M(T) \)
denote the \( \ring  \)-module defined recursively by the rules
\begin{enumerate}
\item if \( T \) is a single leaf-node labeled by a \( \ring  \)-module
\( A \), then \( M(T)=A \).
\item if the root of \( T \) is labeled with \( \homa  \) and its two
children are expression-trees \( T_{1} \) and \( T_{2} \), respectively,
then\[
M(T)=\homz (M(T_{1}),M(T_{2}))\]

\item if the root of \( T \) is labeled with \( \otimes  \) and its children
are expression-trees \( T_{1},\dots ,T_{n} \) then\[
M(T)=\bigotimes _{i=1}^{n}M(T_{i})\]

\end{enumerate}
\end{defn}
\begin{rem}
This associates a multilinear functor of the leaf-nodes with an expression
tree.

For instance, if \( T \) is the expression tree in remark~\ref{rem:expressionstree1},
then \[
M(T)=\homz (A_{1}\otimes A_{2}\otimes A_{3},A_{4}\otimes A_{5})\otimes A_{6}\]

In other words, \( T \) is nothing but the \emph{syntax tree} of
the functors that make up \( M(T). \)
\end{rem}
Now we define \emph{operations} that can be performed on expression
trees and their effect on the associated functors.

Throughout this discussion, \( T \) is some fixed expression tree. 

\begin{defn}
\emph{\label{def:type0transform}Type-0 transformations.} Perform
the following operations or their inverses:
\begin{description}
\item [\( \homa  \)-transform]Given any subtree, \( A \) of \( T \),
replace it by the subtree \begin{figure}[H]\begin{center}\psfrag{Z}{$\ring$}\includegraphics{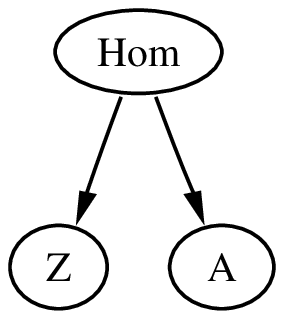} \end{center}\caption{\label{fig:type0trans}}\end{figure}
\item [\( \otimes  \)-transform]Given a subtree of the form \begin{figure}[H]\begin{center}\psfrag{T1}{$T_1$}\psfrag{Ti}{$T_i$}\psfrag{Tip1}{$T_{i+1}$}\psfrag{Tn}{$T_{n}$}
\includegraphics{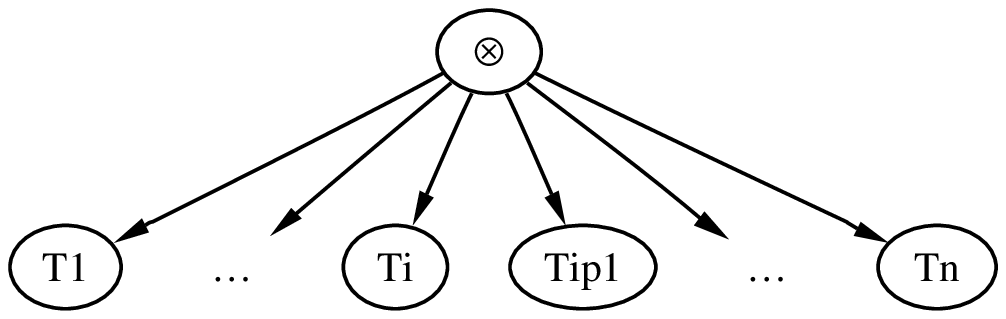} \end{center}\caption{\label{fig:type0timesin}}\end{figure} where
\( n>0 \) is some integer and \( T_{1},\dots ,T_{n} \) are subtrees,
replace it by
\end{description}
\begin{figure}[H]\begin{center}\psfrag{T1}{$T_1$}\psfrag{Ti}{$T_i$}\psfrag{Tip1}{$T_{i+1}$}\psfrag{Z}{$\ring$}\psfrag{Tn}{$T_{n}$}
\includegraphics{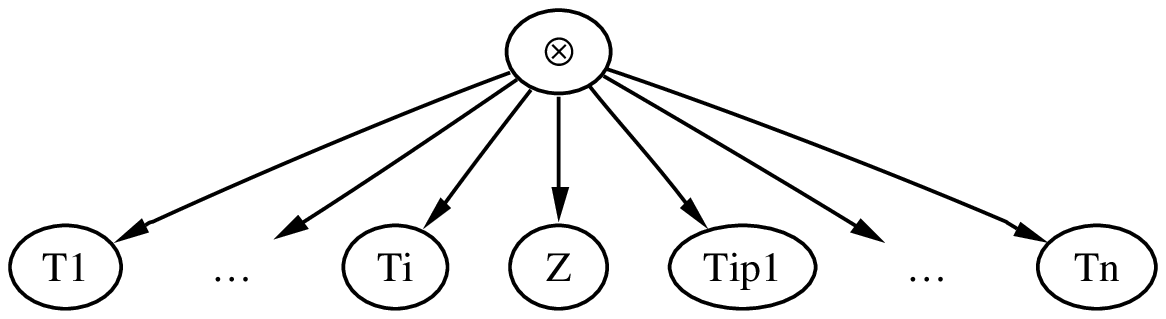} \end{center}\caption{\label{fig:type0timesout}}\end{figure} where
\( 0\le i\le n \)

\end{defn}
In addition, we define slightly more complex transformations

\begin{defn}
\emph{\label{def:type1transform}Type-1 transformations.} Perform
the following operation or its inverse: If \( T \) has a \emph{covariant}
node that is the root of a subtree like \begin{figure}[H]\begin{center}\psfrag{T1}{$T_1$}\psfrag{T2}{$T_2$}\psfrag{T3}{$T_3$}
\includegraphics{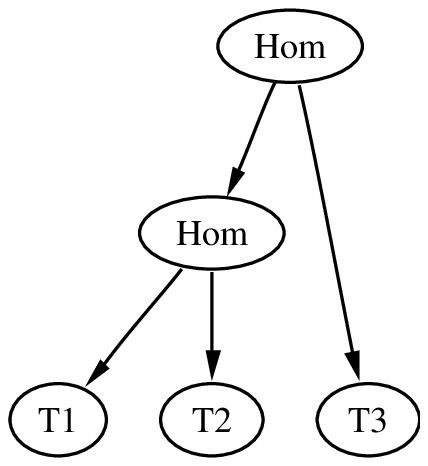} \end{center}\caption{\label{fig:type1covariant}}\end{figure}
where \( T_{1} \), \( T_{2} \), and \( T_{3} \) are subtrees, replace
it by the subtree \begin{figure}[H]\begin{center}\psfrag{T1}{$T_1$}\psfrag{T2}{$T_2$}\psfrag{T3}{$T_3$}
\resizebox*{2in}{2in}{\includegraphics{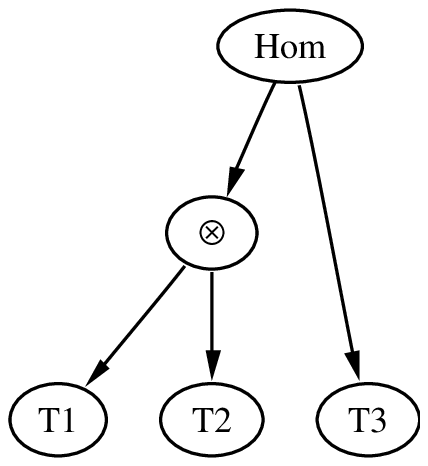}} \end{center}\caption{\label{fig:type1covariantout}}\end{figure}
If it has a \emph{contravariant} node that is the root of a subtree
like figure~\ref{fig:type1covariantout}, replace it by the subtree
depicted in figure~\ref{fig:type1covariant}.
\end{defn}
Finally, we define the most complex transformation of all

\begin{defn}
\emph{\label{def:type2transform}Type-2 transformations.} If \( T \)
is an expression tree with a \emph{covariant} node that is the root
of this subtree like \begin{figure}[H]\begin{center}\psfrag{A1}{$A_1$}\psfrag{An}{$A_n$}\psfrag{B1}{$B_1$}\psfrag{Bn}{$B_n$}
\includegraphics{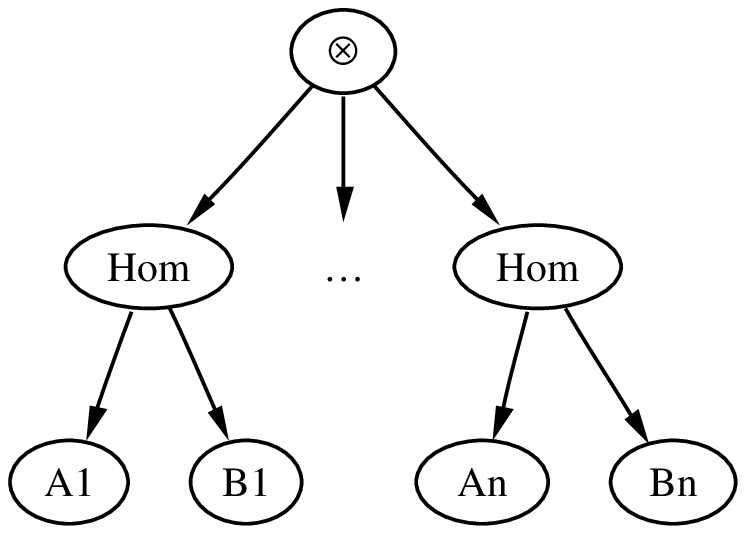} \end{center}\caption{\label{fig:type2covariant}}\end{figure}
where \( n>1 \) is an integer and \( A_{1},\dots ,A_{n} \) and \( B_{1},\dots ,B_{n} \)
are subtrees, we replace the subtree in figure~\ref{fig:type2covariant}
by \begin{figure}[H]\begin{center}\psfrag{A1}{$A_1$}\psfrag{An}{$A_n$}\psfrag{B1}{$B_1$}\psfrag{Bn}{$B_n$}
\includegraphics{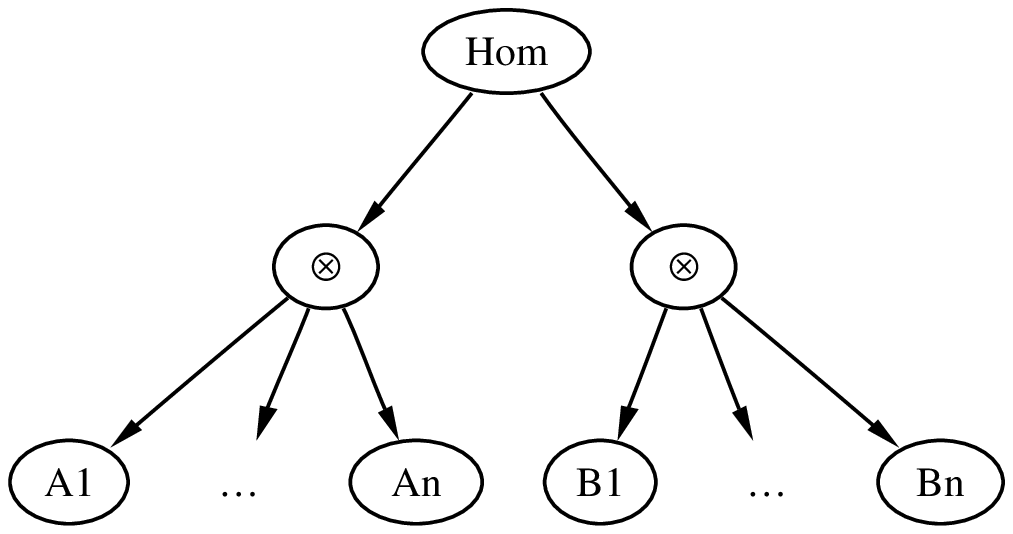} \end{center}\caption{\label{fig:type2covariantout}}\end{figure}
If a node is \emph{contravariant} and is the root of a subtree like
figure~\ref{fig:type2covariantout}, we replace it by the tree in
figure~\ref{fig:type2covariant}.
\end{defn}
Given these rules for transforming an expression tree, we can define
an \emph{induced natural transformation} of functors \( M(T) \):

\begin{claim}
\label{def:exptreeinducedmap}Let \( T \) be an expression tree and
let \( T^{\prime } \) be the result of performing a transform \( e \),
defined above, on \( T \). Then there exists an induced natural transformation
of functors\[
f(e):M(T)\to M(T^{\prime })\]

Given a sequence \( \mathbf{E}=\{e_{1},\dots ,e_{k}\} \) of elementary
transforms, we define \( f(\mathbf{E}) \) to be the composite of
the \( f(e_{i}) \), \( i=1,\dots ,n \).
\end{claim}
This claim follows immediately from the recursive description of \( M(T) \)
in definition~\ref{def:exptreemodule}, the well-known morphisms\begin{eqnarray}
\homz (\ring ,A) & = & A\label{eq:specmorphism1} \\
A\otimes \ring \otimes B & = & A\otimes B\label{eq:specmorpism2} \\
\homz (A,\homz (B,C)) & = & \homz (A\otimes B,C)\label{eq:specmorphism3} \\
\homz (A,B)\otimes \homz (C,D) & \to  & \homz (A\otimes C,B\otimes D)\label{eq:specmorphism4} 
\end{eqnarray}
(where \( A \), \( B \), \( C \), and \( D \) are free \( \ring  \)-modules),
and the functoriality of \( \otimes  \) and \( \homz (\ast ,\ast ) \).

In the case where the \( \ring  \)-modules are DG-modules, we apply
the Koszul convention for type-2 transformations such a transformation
sends\[
(a\mapsto b)\otimes (c\mapsto d)\]
 to\[
(-1)^{\dim b\cdot \dim c}a\otimes c\mapsto b\otimes d\]

The Koszul conventions does \emph{not} produce a change of sign in
any of the other cases.

Now we are ready to state the main result of the appendix:

\begin{thm}
\label{th:linearrigidity}Let \( T \) be an expression tree and suppose
\( \mathbf{E}_{1} \) and \( \mathbf{E}_{2} \) are two sequences
of elementary transformations (as defined in definitions~\ref{def:type0transform}
through \ref{def:type2transform}) that both result in the same transformed
tree, \( T^{\prime } \). Then \[
f(\mathbf{E}_{1})=f(\mathbf{E}_{2}):M(T)\to M(T^{\prime })\]
This result remains true if the \( \ring  \)-free modules on the
leaves are DG-modules and we follow the Koszul Convention.
\end{thm}
\begin{rem}
{}``Same'' in this context means {}``isomorphic.'' This theorem
shows that the induced natural transformation, \( f(\mathbf{E}) \),
\emph{only} depends on the structure of the resulting tree, not on
the sequence of transforms used. There is \emph{less structure} to
maps of the form \( f(\mathbf{E}) \) than one might think.
\end{rem}
We devote the rest of this section to proving this result. We begin
with

\begin{defn}
\label{def:thattransform}Let \( T \) be an expression tree. Then
\( \mathrm{inorder}(T) \) denote the list of leaf-nodes of \( T \)
as encountered in an in-order traversal of \( T \), i.e. 
\begin{enumerate}
\item if \( T \) is a single node \( A \), then \( \mathrm{inorder}(T)=\{A\} \)
\item if the root of \( T \) has child-subtrees \( A_{1},\dots ,A_{n} \)
then\[
\mathrm{inorder}(T)=\mathrm{inorder}(A_{1})\bullet \cdots \bullet \mathrm{inorder}(A_{n})\]
where \( \bullet  \) denotes concatenation of lists.
\end{enumerate}
\end{defn}
Given transformations and in-order traversals, we want to record the
effect of the transformations on these ordered lists.

\begin{prop}
\label{prop:mtlementdesc}Let \( T \) be an expression tree and suppose
the \( \ring  \)-free modules on its leaves are equipped with \( \ring  \)-bases.
Then an element \( x\in M(T) \) can be described as a set of lists\[
x=\{(a_{1},\dots ,a_{k})\dots \}\]
where \( a_{i}\in A_{i} \) and \( A_{i} \) is the free \( \ring  \)-module
occurring in the \( i^{\mathrm{th}} \) node in \( \mathrm{inorder}(T) \).
\end{prop}
\begin{rem}
To actually \emph{define} \( M(T) \) as a free \( \ring  \)-module,
we must add quantifiers and relations that depend on the internal
structure of \( T \) to these lists.
\end{rem}
\begin{proof}
Let \( A \) and \( B \) be free \( \ring  \)-modules. Elements
of \( A\otimes B \) can be described as \( a\otimes b \), where
\( a\in A \), \( b\in B \) are basis elements. So the list in this
case has two elements and the set of lists contains a single element:\[
\{(a,b)\}\]

Elements of \( \homz (A,B) \) are functions from \( A \) to \( B \)
--- i.e., a set of ordered pairs\[
\{(a_{1},b_{1}),\dots ,(a_{i},b_{i},\dots )\}\]
where \( a_{i}\in A \) is a basis element, \( b_{i}\in B \) (not
necessarily a basis element) and \emph{every} basis element of \( A \)
occurs as the left member of some ordered pair. The general statement
follows from the recursive definition of \( M(T) \) in definition~\ref{def:exptreemodule}
and the definition of in-order traversal in definition~\ref{def:thattransform}.

Now we prove theorem~\ref{th:linearrigidity}:

Let \( x\in M(T) \) be given by \[
x=\{(a_{1},\dots ,a_{k})\dots \}\]
as in proposition~\ref{prop:mtlementdesc}. We consider the effect
of the transformations defined in definitions~\ref{def:type0transform}
through \ref{def:type2transform} on this element.
\begin{description}
\item [Type-0]transformations insert or remove terms equal to \( 1\in \ring  \)
into each list in the set.
\item [Type-1]transformations have \emph{no effect} on the lists (they only
affect the \emph{predicates} used to define the module whose elements
the lists represent).
\item [Type-2]transformations permute portions of each list in \( x \).
In the DG case, whenever an element \( a \) is permuted past an element
\( b \), the list is multiplied by \( (-1)^{\dim a\cdot \dim b} \).
\end{description}
Note that, in \emph{no} case is the \emph{data} in the lists altered.
Furthermore, we claim that the equality of the trees resulting from
performing \( \mathbf{E}_{1} \) and \( \mathbf{E}_{2} \) on \( T \)
implies that:
\begin{itemize}
\item the permutations of the lists from the type-2 transformations must
be compatible
\item the copies of \( \ring  \) inserted or removed by the type-0 transformations
must be in compatible locations on the tree.
\end{itemize}
Consequently, the lists that result from performing \( \mathbf{E}_{1} \)
and \( \mathbf{E}_{2} \) on the lists of \( x \) must be the same
and\[
f(\mathbf{E}_{1})(x)=f(\mathbf{E}_{2})(x)\]
 The isomorphism of final expression trees also implies that the predicates
that apply to corresponding element of these lists are also the same.
Since this is true for an \emph{arbitrary} \( x \) we conclude that\[
f(\mathbf{E}_{1})=f(\mathbf{E}_{2})\]

In the \emph{DG case,} we note that type-2 transformations may introduce
a change of sign. Nevertheless, the fact that the elements in the
lists are in the same order implies that they have been permuted in
the same way --- and therefore have the same sign-factor.
\end{proof}
We can generalize (relativize) theorem~\ref{th:linearrigidity} slightly.
We get a result like theorem~\ref{th:linearrigidity} except that
we have introduced a morphism that is \emph{not} of the type

\begin{thm}
\label{th:rigidity2}Let \( T \) be an expression tree whose leaf-modules
are \( \{A_{1},\dots ,A_{n}\} \) and consider the diagram \[\xymatrix{ {} & {T_1}\ar[r]^{\varphi} & {T_3}\ar[rd]^{\mathbf{E}_3} & {}\\ {T}\ar[ru]^{\mathbf{E}_1}\ar[rd]_{\mathbf{E}_2} & {} & {} & {T'} \\ {} & {T_2}\ar[r]_{\varphi} & {T_4}\ar[ru]_{\mathbf{E}_4} & {}}\]
where 
\begin{enumerate}
\item for some fixed index \( k \), \( \varphi :A_{k}\to \bar{T} \) replaces
the leaf node labeled with the module, \( A_{k} \), with an expression
tree \( \bar{T} \) that has leaf-modules \( \{B_{1},\dots ,B_{t}\} \), 
\item \( \mathbf{E}_{1} \), \( \mathbf{E}_{2} \), \( \mathbf{E}_{3} \),
and \( \mathbf{E}_{4} \) are sequences of elementary transformations
(as defined in definitions~\ref{def:type0transform} through \ref{def:type2transform}).
\item \( f(\varphi ):A_{k}\to M(\bar{T}) \) is some morphism of free \( \ring  \)-modules
\end{enumerate}
Then\[
f(\mathbf{E}_{3})\circ f(\varphi )\circ f(\mathbf{E}_{1})=f(\mathbf{E}_{4})\circ f(\varphi )\circ f(\mathbf{E}_{2})\]

\end{thm}
\begin{rem}
In other words, theorem~\ref{th:linearrigidity} is still true if
we have a morphism in the mix that is \emph{not} of the canonical
type in equation~\ref{eq:specmorphism1} through \ref{eq:specmorphism4}
--- as long as the remaining transformations are done in a compatible
way. 
\end{rem}
\begin{proof}
Let \( x\in M(T) \) be given by \[
x=\{(a_{1},\dots ,a_{n})\dots \}\]
We get \begin{eqnarray*}
f(\mathbf{E}_{1})(x) & = & \{(a_{\sigma (1)},\dots ,a_{\sigma (n)})\dots \}\\
f(\mathbf{E}_{2})(x) & = & \{(a_{\tau (1)},\dots ,a_{\tau (n)})\dots \}
\end{eqnarray*}
where \( \sigma ,\tau \in S_{n} \) are permutations. In each of these
lists, we replace \( a_{k} \) by a set of lists\[
\{(b_{1},\dots ,b_{t})\}\]
representing the value of \( \varphi (a_{k}) \) and apply \( f(\mathbf{E}_{3}) \)
and \( f(\mathbf{E}_{4}) \), respectively --- possibly permuting
the resulting longer lists. As in theorem~\ref{th:linearrigidity},
the result is two copies of the same set of lists. This is because
both sets of operations result in the expression tree \( T^{\prime } \),
implying that the permutations must be compatible. As in theorem~\ref{th:linearrigidity},
the key fact is that the data in the lists is not changed (except
for being permuted).
\end{proof}

\providecommand{\bysame}{\leavevmode\hbox to3em{\hrulefill}\thinspace}
\providecommand{\MR}{\relax\ifhmode\unskip\space\fi MR }
\providecommand{\MRhref}[2]{%
  \href{http://www.ams.org/mathscinet-getitem?mr=#1}{#2}
}
\providecommand{\href}[2]{#2}


    \end{document}